\documentclass[12pt]{article}
\usepackage{MnSymbol}
\usepackage{amsmath,amsfonts,
verbatim}

\newcommand{\A}{\mathbb{A}}
\newcommand{\Q}{{\mathbb Q}}
\newcommand{\C}{{\mathbb C}}
\newcommand{\N}{{\mathbb N}}
\newcommand{\Z}{{\mathbb Z}}
\newcommand{\R}{{\mathbb R}}

\newcommand{\HH}{\mathbb{H}}

\newcommand{\F}{\mathrm{F}}

\newcommand{\X}{{\mathbb X}}
\newcommand{\Y}{\mathbb{Y}}

\newcommand{\pr}{{\rm pr}}

\newcommand{\la}{\langle}
\newcommand{\ra}{\rangle}
\newtheorem{pkt}{}[section]  
\newcommand{\bpk}{\begin{pkt}\rm }  
\newcommand{\epk}{\end{pkt}} 
\newcommand{\inv}{^{-1}}   

\newcommand{\be}{\begin{equation}}  
\newcommand{\ee}{\end{equation}}

\newcommand{\Ss}{\mathrm{S}}
\newcommand{\kk}{\mathrm{k}}
\newcommand{\G}{\mathrm{G}}
\newcommand{\GL}{\mathrm{GL}}
\newcommand{\SL}{\mathrm{SL}}

\newcommand{\Aut}{\mathrm{Aut}}
\newcommand{\Gal}{\mathrm{Gal}}

\newcommand{\ttt}{\mathbf{t}}
\newcommand{\rrr}{\mathbf{d}}

\usepackage{graphicx}
\usepackage[all,cmtip]{xy}
\usepackage{tikz}

\title{Canonical models of modular curves and the Galois action on CM-points}
\author{C.Daw\footnote{Supported by EPSRC New Investigator Award EP/S029613/1} and B.Zilber\footnote{Supported by EPSRC Program Grant ``Symmetries and Correspondences''}  }
 
\begin{document}

\maketitle

We use the theory of canonical models of Shimura varieties to describe the projective limit of the curves $\Y(N),$ all $N,$ and its automorphism group. 
In particular we prove that the Galois group of $\Q$(CM) over $\Q$ is an open subgroup of 
an automorphism group of an explicitly given adelic structure.

\section{Introduction}
The origin of this paper is due to a model-theoretic project which attempts to formalise the notion of an analytic covering space of an algebraic variety in such a way that the formal cover is unique up to isomorphism. See e.g. \cite{DH} or \cite{Zspecial} for a survey on this project. In particular, we are interested in \cite{CZ} in the situation that arises in the context of a system of analytic functions
\be\label{j} j_N:\HH \to \Y(N),\ N=1,2\ldots\ee
where $\Y(N)=\Gamma(N)\backslash \HH$ is the algebraic curve isomorphic to the modular curve $\Y(N)$ but defined over $\Q,$ at the cost of ignoring the 
level $N$ structure on the modular curve. Such curves (or rather their completed versions $\X(N)$) are considered e.g. in \cite{Xavier}.

The key to obtaining the desired uniqueness result turns out to be the understanding of $$\tilde{\HH}:= \varprojlim \Gamma(N)\backslash \HH$$
and the action of $\Aut\,\C$ on it. The latter is reduced to the understanding of
$\Gal(\Q(\mathrm{CM}): \Q),$ where CM is the set of all the CM points on all the $\Y(N).$ More precisely, we need to describe the
structure, the algebraic relations on the set of CM-points and the action of  $\Gal(\Q(\mathrm{CM}): \Q)$ on CM-points.

\medskip     

We use the well-known machinery of canonical models of Shimura varieties and complex multiplication, see e.g. \cite{Milne0}, which is essentially an answer to the Langlands conjecture on conjugation of Shimura varieties in \cite{Langlands}.
\bpk
A characterisation of $\tilde{\HH}$ is closely related to a problem raised in \cite{Langlands} by R.Langlands that was resolved by M.Borovoi, P.Deligne, J.Milne, K.Shih and others. (We only need the case of the product of Shimura curves.)
A key outcome of the theory is the proof of existence and uniqueness of {\bf canonical models} of Shimura varieties.

Recall that a Shimura variety is given by a Shimura datum $(G,X)$, where $X$ is a complex manifold and $G$ is a reductive linear algebraic group such that $G(\R)$ acts on $X.$
The corresponding Shimura variety is defined as
$$\mathrm{Sh}(G,X)_K:=G(\Q)\backslash(X\times G(\A_f)/K)$$
where $\A_f$ is the ring of finite adeles and $K$ runs through the family of open compact subgroups of $G(\A_f).$ For each $K,$ the set $\mathrm{Sh}(G,X)_K$ of double cosets naturally possesses the structure of a complex quasi-projective algebraic variety.  The canonical model, which we also denote $\mathrm{Sh}(G,X)_K$, is a model for $\mathrm{Sh}(G,X)_K$ over an explicit number field $E(G,X)$. For $X=\HH^{\pm}:= \C\setminus \R$ and $G=\GL_2$ 
we get the family of modular curves $\mathrm{Sh}(G,X)_K$ defined over $\mathbb{Q}$. We refer to \cite{Milne0} for a general survey and \cite{Milne1} for the case of Shimura curves. 

 Denote
$$\Ss:=\GL_2(\Q)\backslash (\HH^{\pm}\times \GL_2(\A_f)),$$
{\em the covering of canonical Shimura curves},
and $\Delta$ the subgroup of matrices of $\GL_2(\A_f)$ of the form $\left(\begin{array}{ll}
a\ 0\\ 0\ 1
\end{array}\right),$ $a\in \hat{\Z}^\times.$ Then there is a holomorphic covering map $\Ss\to \tilde{\HH}$ which represents $\tilde{\HH}$ as  $$\tilde{\HH}\cong \Ss/\Delta= \GL_2(\Q)\backslash (\HH^{\pm}\times \GL_2(\A_f))/\Delta.$$
This essentially reduces the study of $\tilde{\HH}$ to the study of
$\Ss.$
\epk
\bpk In our analysis it is important to distinguish the full structure on $\Ss$ as a complex manifold with infinitely many connected components, which we name $\mathbf{S}^\mathrm{Full},$ and the purely algebraic structure  $\mathbf{S}^\mathrm{Pure},$ which is just a set split into connected components with the action of $\GL_2(\A_f)$ on it. 

It is clear that a CM-point on the $\Y(N)$ comes from a point on $\Ss$  fixed by a $g\in \GL_2(\A_f)$ which is not in the centre of $\GL_2(\Q)$ (the centre  of $\GL_2(\Q)$ acts trivially on $\Ss$). We call such points fixed points or, with a slight abuse of terminology, CM-points. 

In particular we are interested in substructures 
$\mathbf{S}^\mathrm{Full}(\mathrm{CM})$ and $\mathbf{S}^\mathrm{Pure}(\mathrm{CM})$ of our structures reduced to the sets of their CM-points (also sometimes called {\em special points}). 

\medskip

Note that by construction $\Ss$ is a union of $\GL_2(\A_f)$-orbits (Hecke orbits) $S_\tau:=\GL_2(\A_f)\cdot \tau,$ $\tau\in \HH^{\pm}$ 
 and in particular 
$$\mathrm{CM}=\bigcup_{\sqrt{-m}}\mathrm{CM}_{\sqrt{-m}},\ \ \mathrm{CM}_{\sqrt{-m}}=\GL_2(\A_f)\cdot \sqrt{-m}$$
where $m$ runs in $\N,$ square-free (including $m=1$). 
\epk

\bpk The main results are as follows. 

{\bf Theorem A} (see \ref{MainGal}(i)). $$\Gal(\Q(\mathrm{CM}):\Q)\cong\Aut\, \mathbf{S}^\mathrm{Full}(\mathrm{CM})\le_\approx \Aut\, \mathbf{S}^\mathrm{Pure}(\mathrm{CM}),$$
where $\le_\approx$ ("almost equal")
 means that for any finite set  $\bar{\tau}$ of square-free imaginary quadratics, for $\mathrm{CM}_{\bar{\tau}}=\bigcup \{ \mathrm{CM}_\tau: \tau\in \bar{\tau}\},$
$$ \Aut\, \mathbf{S}^\mathrm{Full}(\mathrm{CM}_{\bar{\tau}})\le_\mathrm{open} \Aut\, \mathbf{S}^\mathrm{Pure}(\mathrm{CM}_{\bar{\tau}}),$$ 
(open, and hence finite index subgroup).

The theorem is the basis of the model-theoretic result in our paper \cite{CZ} which we expect to be extendable to the general case of Shimura varieties, and so the statement of Theorem A.

The  theorem allows, with some extra work, see \ref{Gapprox},  to characterise the limit structure $\tilde{\HH}=\varprojlim \Gamma(N)\backslash \HH$ and the abstract formalisation of the system (\ref{j}), which is the key  to the application in \cite{CZ}. 

\medskip   
Let $\F$ be
the composite of all $\Q(\tau)^\mathrm{ab}$ for all imaginary quadratic  $\tau.$
Let $P\subset \N$ be the set of primes, including 1. 

\medskip

{\bf Theorem B} (see \ref{MainGal}(ii)).

 (i) $$\Gal(\Q(\mathrm{CM}):\Q)\cong \Gal(\F:\Q)$$

(ii)
$$\Gal(\F:\Q^\mathrm{ab})\le _\approx \prod_{p\in P} \Gal(\Q(\sqrt{-p})^\mathrm{ab}: \Q^\mathrm{ab})$$
(the RHS can be identified with a subgroup of $\Aut\, \mathbf{S}^\mathrm{Pure}(\mathrm{CM})$) and  $$\Gal(\Q(\sqrt{-p})^\mathrm{ab}: \Q^\mathrm{ab}) \cong T^0_{\sqrt{-p}}(\A_f)/(T^0_{\sqrt{-p}}(\A_f)\cap \GL_2(\Q))$$ for some 1-dimensional torus  $T^0_{\sqrt{-p}}$ over $\Q.$

\medskip

(iii) There  is a commutative group $\mathbb{T}$ given in terms of an infinite dimensional matrix group over finite adeles $\A_f,$
$$\mathbb{T}\le_\approx \prod_{p\in P} \Gal(\Q(\sqrt{-p})^\mathrm{ab}: \Q(\sqrt{-p}))$$
and a short exact sequence
$$1\to \mathbb{T}\to \Gal(\Q(\mathrm{CM}):\Q)\to \mathrm{C}_2^P\to 1$$
where $\mathrm{C}_2^P$ is the infinite product of 2-element groups $\mathrm{C}_2(p),$ $p\in P,$ each acting by an involution on the $p$-th co-ordinate  of $\mathbb{T}.$

\medskip

We are especially interested in $\tilde{\HH}$ with the structure induced on it from $\mathbf{S}^\mathrm{Full},$ including the adelic group $\tilde{\G}< \GL_2(\A_f)$ acting on $\tilde{\HH}$ by the action which comes from the action of $\GL_2(\A_f)$ on $\Ss.$
This group can 
 be seen as the completion of $\GL_2^+(\Q)$ in $\tilde{\HH}.$ In particular, we establish that $\Aut\, \C$ acts on $\tilde{\G}$ as the Galois group  $\Gal(\Q^\mathrm{ab}:\Q)\cong \hat{\Z}^\times,$ by conjugation with elements of $\Delta.$

\epk

\bpk \label{Intro1.5} Note that 
Theorem A brings into the subject a Model Theory flavour by focussing on the specific class of structures $\mathbf{S}^\mathrm{Pure}$ which in Model Theory classification belong to the class of {\em locally modular geometries of trivial type}. It thus emphasises the analogy with Serre's open image theorem which characterises 
$\Gal(\kk(\mathrm{Tors}):\kk)$ of torsion points of an elliptic curve $E$ (without CM) defined over $\kk$ as an open subgroup in the linear group $\GL_2(\hat{\Z}).$ That is it characterises the structure $E(\mathrm{Tors})$ as a {\em locally modular geometry of linear type, an analogue of $\mathbf{S}^\mathrm{Pure}$}.

\epk
\bpk 
The first author would like to thank the EPSRC for its support via a New Investigator Award (EP/S029613/1). He would also like to thank the University of Oxford for having him as a Visiting Research Fellow.

The second author thanks EPSRC for support provided through the Program Grant "Symmetries and Correspondences" and also wishes to express his gratitude to J.Derakhshan for his interest in the work and his help with some of the mathematical issues in the paper.

\epk 
\section{Canonical models of Shimura varieties}
\bpk
We follow \cite{Milne0}, Ch.12, and \cite{Milne1}, section 2. The latter has more specific  details in the case of modular curves although we also need some theory of canonical models of  Shimura varieties of a more general form.

Recall that a Shimura variety $\mathrm{Sh}(G,X)$ is given by the Shimura datum $(G,X)$ where $X$ is a complex manifold and $G$ is a reductive linear algebraic group defined over $\Q$ together with an action of $G(\R)$  on $X$ which satisfies certain conditions (SV1--SV3, see \cite{Milne0}, Definition 5.5. In the full generality of the conjecture of \ref{Intro1.5} we also assume SV5). 

Our main interest is in the case when $$G=\GL_2\mbox{ and }
X=\HH\cup -\HH=\C\setminus \R$$
but we also  need the cartesian powers of this Shimura datum,
$(G^n, X^n)$, which are again Shimura data by definition.
\epk
\bpk \label{introSh}
Given a compact open subgroup $K$ of $G(\A_f)$
one defines the canonical  model for respective Shimura datum to be the scheme of an appropriate (reflex) number field 
$$\mathrm{Sh}(G,X)_K:=G(\Q)\backslash(X\times G(\A_f))/K$$
which should be understood as the set of double cosets with scheme structure (which we use in a very restricted way) induced  by the ring of automorphic forms on $X.$ This definition also makes sense when $K$ is not open, e.g. $K=\{ 1\}.$ In this case it is still a scheme though not of finite type. See \cite{Milne0}, Remark 5.30. We are going to use the construction for some compact $K$ which are not open in this more general sense without making claims on its nature as a scheme.  

When $K$ is an open compact in $G(\A_f)$ then $\mathrm{Sh}(G,X)_K$ is a complex algebraic variety.

When $K=\{ 1\}$ we get the complex manifold  $$\mathrm{Sh}(G,X):=G(\Q)\backslash(X\times G(\A_f))=\lim_K \mathrm{Sh}(G,X)_K $$
 with uncountably many connected components (this equality assumes certain assumptions which $G=\GL_2$ satisfies).

\medskip

The elements  $(x,a)\in X\times G(\A_f)$ are acted upon by $q\in G(\Q)$ on the left and by $k\in K$ on the right, so 
$$q\cdot (x,a)\cdot k:=(qx,qak),$$ and the respective elements of $\mathrm{Sh}(G,X)_K$ are written in the form $[x,a]$ and by definition of double cosets
 $$[x,a]= [qx,qak].$$

It will be convenient for us to define
 the action of $g\in G(\A_f)$ on the set of cosets, for $K=\{ 1\},$ by multiplication by $g\inv$ on the right
\be\label{gs} g* [x,a]:= [x,ag\inv].\ee
If $g=q\in G(\Q)$ then we have
$$ q* [x,1]=[x,q\inv]=[q x,1],$$
 that is, if we identify $X=X\times \{ 1\}\hookrightarrow \mathrm{Sh}(G,X)_K$  with $K=\{ 1\}$ then
 $G(\Q)$ acts on $X$ exactly by its standard action. 
 In particular, 
 \be \label{qs}  q x=x\Rightarrow  q* [x,1]=[x,1].\ee

\epk

\bpk We now fix our notations. For a modular curve we have
\be\label{SK} \Ss_K:= \GL_2(\Q)\backslash  (\HH\cup -\HH) \times \GL_2(\A_f)/K,\ee
where  $K$ is an open compact subgroup of  $\GL_2(\A_f).$ 
For  the case  $K=\{ 1\},$ 
\be\label{S} \Ss:= \GL_2(\Q)\backslash  (\HH\cup -\HH) \times \GL_2(\A_f).\ee

Here $q\in \G(\Q)$ acts on pairs $\la \tau,a\ra\in (\HH\cup -\HH) \times \GL_2(\A_f)$ from the left diagonally,  $\la \tau,a\ra\mapsto \la q\tau,qa\ra,$  and $k\in K$ acts on the right on the second coordinate, $\la \tau,a\ra \mapsto \la \tau,ak\ra.$ Thus an element of $\Ss_K$ is a class $[\tau,a]_K$ defined by the condition 
\be\label{qak} [\tau,a]_K=[\tau',a']_K \Leftrightarrow (q\tau, qak)=(\tau',a'),\mbox{ for some }q\in \GL_2(\Q),\ k\in K .\ee

\epk
\bpk\label{Del} {\bf Remark.} Note that according to Deligne \cite{Deligne1979}, 
$\Ss_K,$ for $K$ compact open in $\GL_2(\A_f)$  are complex algebraic curves defined over $\Q$ (usually geometrically reducible). 
We are specifically interested in subgroups
$$K(N):=\Delta(N)\cdot \tilde{\Gamma}(N), \mbox{ where } \Delta(N)= \{ \left(\begin{array}{ll}
a\ 0\\ 0\ 1
\end{array}\right): \ a\in \hat{\Z}^\times,\ a\equiv 1 \mod N\}$$
and $$\tilde{\Gamma}(N)=\{ \left(\begin{array}{ll}
a\ b\\ c\ d
\end{array}\right)\in  \SL_2(\hat{\Z}): \ a\equiv d\equiv 1\,\mbox{mod}\, N,\ b\equiv c\equiv 0\,\mbox{mod}\, N\}.$$
 Each component of $\Ss_{K(N)}$  together with the level $N$ structure is defined over $\Q(\zeta_N),$
 and is isomorphic to
$\Y(N)_{/\Q(\zeta_N)},$ the level $N$ modular curve over $\Q(\zeta_N).$
 Here and below by ``the level $N$ structure" we understand the action of individual elements $\gamma$ of $\Gamma$ on $\Y(N).$ This defines the curve $$C_{\gamma,N}=\{ \la j_N(\tau), j_N(\gamma \tau)\ra: \tau\in \HH\} \subset \Y(N)\times \Y(N)$$
which is definable over $\Q(\zeta_N).$  
(We later apply the notation $C_{g,N}$ to define an algebraic curve in $\Y(N)\times \Y(N)$ corresponding to the general $g\in \GL_2(\A_f)$).

 Moreover, the components are conjugated by Galois automorphisms over $\Q$ and the whole (geometrically reducible) curve is defined over $\Q.$ 
See more on this in \ref{defApprox}.

 We will write $\Delta$ for $\Delta(1)$ below.
When $K=K(N)$ the structure $\Ss_{K(N)}$  is an irreducible algebraic curve over $\Q$ and isomorphic to $\Y(N).$ Taking the projective limit over $N$ we get  
$ \Ss_{\Delta}$ (also denoted $\Ss_\approx$ later) which is isomorphic by construction to $\tilde{\HH}.$ 

\epk

\bpk\label{definitionS} {\bf Definition.}
 We represent the family $\{ \mathrm{S}_K\}$  as  the multisorted structure
  $\mathbf{S}$ with sorts $\Ss$ and  $\Ss_{K(N)},$ $N\in \N,$ such that:
 \begin{itemize}
 \item[(i)] $\mathrm{S}$ is the set with unary operations   $s\mapsto g* s,$  for each  $g\in\GL_2(\A_f);$ 
 \item[(ii)]  the surjective maps: $\pr_K: \mathrm{S}\to \mathrm{S}_K$ are definable in $\mathbf{S},$ and so are the maps $$\pr_{K,K'}: \mathrm{S}_K\to \mathrm{S}_{K'}\mbox{ for } K=K(N),\ K'=K(N'),\ N'|N. $$
  \item[(iii)] the partition of algebraic curves $\Ss_K$ into irreducible components and the respective partition (equivalence relations $E_K$ and $E$) of $\Ss$ are definable. 
 \item[(iv)] the relations and operations induced from (i)-(iii)  on $\mathrm{S}_K$ are definable.  
    
\end{itemize}     

Note that we ignored the complex structure on the $\mathrm{S}_K.$ We call the so defined structure $\mathbf{S}$ (or $\mathbf{S}^\mathrm{Pure}$) the {\bf pure levels structure}.
The {\bf pure level $N$ structure} is $\mathrm{S}_{K(N)}.$ 

The pure levels structure is part of the {\bf full structure} $\mathbf{S}^\mathrm{Full}$ which includes all  relations and points defined over $\Q$ of the complex algebraic curves $\mathrm{S}_K.$    
\epk

\bpk \label{Fact2.6} {\bf Proposition.} (i) {\em For $g\in \GL_2(\hat{\Z}),$
the operation $s\mapsto g*s$ in \ref{definitionS}(i) is a morphism over $\Q.$ }

(ii) {\em The morphisms $\pr_K$ and $\pr_{K,K'}$ in \ref{definitionS}(ii) are definable over $\Q.$}

(iii) {\em The curves $\Ss_{K},$ for $K=K(N),$ $N\in \N,$ and $K=\{ 1\}$ are defined over $\Q$ and so are the equivalence relations 

$E_K(s_1,s_2) \Leftrightarrow s_1,s_2 \mbox{ belong to the same irreducible component of } \Ss_K.$ }  

\smallskip

{\bf Proof.} We use \cite{Milne0} again.
Note first that canonical models for $\mathrm{Sh}_{K(N)}(\GL_2, \HH)$ over the field $E(\GL_2, \HH)=\Q$ exist according to Deligne \cite{Deligne}, see also \cite{Milne1}. In particular,  $\Ss_{K(N)}$ are defined over $\Q.$ (iii) follows sinse $\Ss$ is the limit of the $\Ss_{K(N)}$ and $E_K$ just defines the splitting into irreducible components of a curve over $\Q.$

  For $g\in \GL_2(\hat{\Z})$ we have 
$g\cdot K(N)\cdot g\inv=K(N)$ and hence by \cite{Milne0}, Theorem 13.6,
the maps 
$$s/K(N)\mapsto g*s/K(N)$$
are defined over $\Q.$  Since the action of $g$ on $S$ is the limit of its actions on the $S_{K(N)}$
 (i) follows.

(ii)  follows from  \cite{Milne0}, Theorem 13.6 when one considers $g=1.$

   $\Box$
   
   \medskip
   
   { Commentary.} In fact, the operation $s\mapsto g*s$ on $S$  is definable over $\Q$ for any $g\in \GL_2(\A_f).$ This is proved in the next section.
\epk

\bpk \label{2.5+} {\bf Remark.} The parts (ii) and (iv) of
 the definition \ref{definitionS} of $\mathbf{S}^\mathrm{Pure}$ can be dropped when we assume that any quotient-set and relation interpretable in the structure (i.e. invariant under automorphisms of the structure) is definable in the structure. 
 
 Indeed, the $\Ss_K$ are such quotient-sets and $\pr_K$ and $\pr_{K,K'}$ can be seen as relations between $\Ss, \Ss_K$ and $\Ss_{K'}.$ 
 
 In (iii) it is enough to have the equuivalence relation $E$ on $\Ss$ since $E_K$ is the image of $E$ under the quotient-map $\Ss\to \Ss_K.$     
\epk
\section{ $\Aut\, \C$ action on $\Ss.$}

\bpk\label{defApprox}
Below, $K=K(N)$ if not stated otherwise.
 
Irreducible components of the algebraic curve $\mathrm{S}_{K}$  correspond (see \cite{Milne1}, 2.5) to  cosets $\SL_2(\A_f)\cdot b\subset \GL_2(\A_f),$ 
which in their own turn  are in a bijective correspondence with  elements of $\Q^*\backslash (\{ \pm 1\}\times \A_f^\times)/\det K(N)$  (\cite{Milne1}, 2.6), where 
$$\det K(N)=\{ \det g: g\in K(N)\}= \{ \mu\in \hat{\Z}^\times: \mu\equiv 1\mod N\}.$$
An irreducible component we  write  respectively as   

\be\label{SKir} \mathrm{S}^{\det{b}}_K:= \mathrm{SL}_2(\Q)\backslash  (\HH \times \SL_2(\A_f)\cdot b)/(K\cap \SL_2(\A_f))\ee

In particular, for $K=\{ 1 \},$ irreducible components of $\mathrm{S}$ are in bijective correspondence with   $$\Q^*\backslash \{ \pm 1\}\times \A_f^\times\cong \hat{\Z}^*$$
and $b$ can be chosen to be 
$$b=\rrr_\lambda=( \begin{array}{ll}
\lambda\ 0\\
0\ 1
\end{array}),\ \ \lambda\in \hat{\Z}^*.$$

Now we define a bijection between components  $\mathrm{S}^\mu,$ $\mu=\det b,$ of $\mathrm{S}:$

$$\phi_\lambda: [\tau,a]\mapsto \rrr_\lambda*[\tau,a]= [\tau, a\rrr_\lambda\inv ],\ \ \mathrm{S}^{\mu}\to \mathrm{S}^{\mu\lambda},$$
$\lambda\in \hat{\Z}^*.$ Note that $\rrr_\lambda$ normalises $K=K(N)$.
Thus the bijection $\phi_\lambda$ induces the bijection
$$\phi_{\lambda,K}: [\tau,aK]\mapsto [\tau, a\rrr_\lambda\inv K],\ \ \mathrm{S}^{\mu}_K\to \mathrm{S}^{\mu\lambda}_K, \ \ \mathrm{S}_K\to \mathrm{S}_K.$$
By \ref{Fact2.6} $\phi_\lambda$ and $\phi_{\lambda,K}$  are bi-regular morphisms defined over $\Q.$

Now we can also define embeddings of $\HH$ into each of the components of $\Ss:$
\be \label{HinS} \mathbf{i}_\mu: \tau\mapsto [\tau,\rrr_\mu],\ \ \HH\hookrightarrow \Ss^\mu.
\ee 
 Consequently, $\mathbf{i}_\mu$ induces
the embeddings 
$$\Gamma(N)\backslash \HH\to S^\mu_{K(N)}$$
which has to be a bijection preserving the complex structure, so a bi-regular isomorphism.

\epk
 \bpk \label{S0+}
 The transformation group $\{ \phi_\lambda: \lambda\in \hat{\Z}^*\}$
 gives rise to the following equivalence relations $\approx_K$ on the $\mathrm{S}_K,$ for $K=K(N):$
 $$[\tau,aK]\approx_K [\tau,a\rrr_\lambda\inv K]= [\tau,aK\rrr_\lambda\inv].$$
It is clear that two distinct points on the same component $\mathrm{S}^\mu_K$ are not equivalent and a component $\mathrm{S}^\mu$ intersects each equivalence class of a point. Thus
$\mathrm{S}_{\approx_K}$  is an algebraic curve bi-regularly  isomorphic to each of the components.

We conclude that
\be \label{SY} 
\mathrm{S}_{\approx_{K(N)}}\cong\ \Y(N)= \Gamma(N)\backslash\HH\ee
where $\Y(N)$ is the irreducible algebraic curve over $\Q$ mentioned in the introduction.
Respectively, we can identify
\be\label{Sapprox} \Ss_\approx=\lim_{\leftarrow}\Gamma(N)\backslash\HH=:\tilde{\HH}\ee
as the set or as the structure.
\epk
\bpk {\bf Notation.} For $\tau\in \HH$ let
 $$G_\tau:=\{ q\in \GL_2(\Q): q\tau=\tau\},$$ \ \ 
$$T_\tau =\{ z\in \GL_2(\A_f): \forall q\in G_\tau\, zq=qz\}$$  and  $$ N_\tau =\{ w\in \GL_2(\A_f): w G_\tau w\inv=G_\tau\}.$$ 

\epk
 \bpk \label{S+}
 Note also that for $\tau=\sqrt{-n},$ $n>0$ square-free, simple calculations show that
 $$ G_{\sqrt{-n}}=\left\lbrace \left(\begin{array}{ll}
 \ a\ nb\\
 -b\ \ a
\end{array}\right):\  a,b\in \Q,\  a^2+nb^2\neq 0\right\rbrace,$$    
 $$ T_{\sqrt{-n}} =\left\lbrace \left(\begin{array}{ll}
 \ a\ nb\\
 -b\ \ a
\end{array}\right):\ a,b\in \A_f,\  a^2+nb^2\in \A_f^\times\right\rbrace.$$
Hence $$N_{\sqrt{-n}} = \left\lbrace \left(\begin{array}{ll}
 \ a\ \ \ \ nb\\
 \mp b\  \pm a
\end{array}\right): \ a,b\in \A_f,\ a^2+nb^2\in \A_f^\times\right\rbrace=\left\langle T_{\sqrt{-n}} , \rrr_{-1}\right\rangle,$$
$N_{\sqrt{-n}} /T_{\sqrt{-n}} $ is a  2-element group.

Note, $T_\tau =T(\A_f)$ for $(T,\tau)$  a special torus in the terminology of \cite{Milne0}.

We will use the notations $T_\tau$ and $N_\tau$ for the subgroups of the group scheme $\GL_2.$ Then  $T_\tau$ is a 2-dimensional torus over $\Q$  and $N_\tau$ its normaliser in $\GL_2.$ 

\medskip

An arbitrary quadratic $\tau\in \HH\cup -\HH$ can be presented as $\tau=s* \sqrt{-n},$ for some $s\in \GL_2(\Q),$ which also presents $G_\tau,T_\tau$ and $N_\tau$ as conjugates of the groups above.  
\epk
\bpk \label{specpoints} {\bf Special points of\ $\mathrm{S}_{K(N)}$.}

We need first to classify fixed points of elements $g\in \GL_2(\A_f)$ acting on $\Ss.$ Note that if $q\in Z(\Q),$ the centre of $\GL_2(\Q),$ then $q\tau=\tau,$ for any $\tau\in \HH.$

Let $\tau\in \HH$ be a fixed point of $q\in \GL_2(\Q)\setminus Z(\Q).$ By (\ref{qs}) also $q\cdot [\tau,1]=[\tau,1].$

Let $a,g\in \GL_2(\A_f).$
 Consider $[\tau,a]\in \Ss.$ 
 We have  \be \label{fixg} g* [\tau,a]=[\tau,a] \Leftrightarrow g\inv=a\inv q'a,\mbox{ for some } q'\mbox{ such that } q' \tau=\tau\ee

Let $\tau\in \HH,$ $a\in \GL_2(\A_f)$ and $q\in \GL_2(\Q)\setminus Z(\Q).$ We claim: 
\begin{itemize}
\item[(i)] $[\tau,a]$ is a fixed point of some non-central  element $g\in \GL_2(\A_f)$ if and only if $\tau\in \HH$ is quadratic and $g\in a\inv G_\tau a$;


\item[(ii)] if $q\tau=\tau$ then $[\tau,a]$ is a fixed point of $q$ iff there is   $q'\in G_\tau$ such that $a qa\inv= q'$ iff $a\in N_\tau ;$

\item[(iii)] $[\tau,aK]$ are special points on $\mathrm{S}_K$ if and only if  $\tau\in \HH$ is quadratic.

\end{itemize}

The proof of (i)--(iii) is by simple calculations. Let us in particular show the last part of (ii): $q=a\inv q'a \Rightarrow a\in N_\tau.$

Shifting $\tau$ by an element of $\GL_2^+(\Q)$ we may assume $\tau=\sqrt{-m}$ for some square-free positive integer $m.$ Then
$$q=\left( \begin{array}{ll}
\ s\ mt\\ -t\ \ s
\end{array}\right),\ \  q'=\left( \begin{array}{ll}
\ s'\ mt'\\ -t'\ \ s'
\end{array}\right)
$$
and
$$q=a\left( \begin{array}{ll}
s\ mt\\ -t\ s
\end{array}\right)a\inv= \left( \begin{array}{ll}
s\ 0\\ 0\ s
\end{array}\right)+a\left( \begin{array}{ll}
0\ mt\\ -t\ 0
\end{array}\right)a\inv
$$
and thus 
$$a\left( \begin{array}{ll}
\ 0\ mt\\ -t\ 0
\end{array}\right)a\inv =\left( \begin{array}{ll}
s''\ mt'\\ -t'\ s''
\end{array}\right)\mbox{ where } s''=s'-s.
$$
We may assume $t=1.$
Since
$\left( \begin{array}{ll}
0\ m\\ -1\ 0
\end{array}\right)^2=\left( \begin{array}{ll}
-m\ 0\\ 0\ -m
\end{array}\right)$ we have $\left( \begin{array}{ll}
s''\ mt'\\ -t'\ s''
\end{array}\right)^2=\left( \begin{array}{ll}
-m\ 0\\ 0\ -m
\end{array}\right)$ which implies $s''=0$ and $t'=\pm 1.$ The case $t=1$ means that $a$ commutes with $\left( \begin{array}{ll}
\ 0\ m\\ -1\ 0
\end{array}\right)$ which implies $a\in T_\tau.$ In the case $t=-1$ the matrix $\rrr_{-1}a$ commutes with $\left( \begin{array}{ll}
\ 0\ m\\ -1\ 0
\end{array}\right)$ and so $a\in \rrr_{-1}T_\tau=N_\tau.$   
\medskip

We will refer to fixed points on $\Ss$  as CM-points, or special points on $\Ss$ and denote the set of special points as $\Ss(\mathrm{CM})$ and, respectively, $\Ss_K(\mathrm{CM}),$ for an open compact subgroup $K.$ 
The same terminology will be used in cases when $K$ is not open.
\epk

\bpk
 Let $\sigma\in \Aut\, \C,$ a field automorphism. This acts naturally on  the complex curves $\Ss_{K(N)}$ and thus on their limit  $\Ss.$ 
 So
$\sigma$ induces also a transformation of $\Ss$ which we keep refering to as $\sigma.$

Below we will assume that  $\sigma$ preserves the pure levels structure which by definition implies, for each
$[\tau,a]\in \Ss$ and $g\in  \GL_2(\A_f),$
\be\label{autg} (g*[\tau,a])^\sigma= g*[\tau,a]^\sigma.\ee

For further analysis we note that $\Ss$ splits into a disjoint union of {\bf Hecke orbits} $$\mathrm{S}_\tau:=[\tau,\GL_2(\A_f)]\subset \mathrm{S}.$$
\epk
\bpk\label{orbits} {\bf Lemma.} {\em Given a quadratic $\tau\in \HH,$ the Hecke orbit $\mathrm{S}_\tau$ can be represented as $$\Ss_\tau= \Ss_{\sqrt{-n}}$$ for some 
 positive square-free integer $n,$ and this representation is unique.}
 
 {\bf Proof.} By definition $\tau= a\sqrt{-n}+b$ for some  $a,b\in \Q$ and a  positive square-free integer $n.$ Then
 $\tau= r \sqrt{-n}$ for $r=\left( \begin{array}{ll}
 a\ b\\
 0\ 1
\end{array}\right),$ that is $\tau \in  \Ss_{\sqrt{-n}}.$  

In order to prove uniqueness consider the assumption that
$[\sqrt{-m},g_1]=[\sqrt{-n},g_2]$ for a square-free positive integer $m$ and $g_1,g_2\in \GL_2(\A_f).$ This is equivalent to 
$[\sqrt{-m},1]=[\sqrt{-n},g],$ for $g=g_2g_1\inv.$  By (\ref{qak}) the latter is equivalent to 
$$(\sqrt{-m},1)=(q\sqrt{-n},qg),\mbox{ for some }q\in \GL_2(\Q),$$
that is $\sqrt{-m}=q\sqrt{-n}.$ This is possible only if $m=n.$ 

The uniqueness of representation follows. $\Box$

\epk
\bpk {\bf Automorphism group $\Aut\,\mathbf{S}_\tau$ in the pure levels structure  $\mathbf{S}^\mathrm{Pure}.$}

The definition of level structure assumes that $\sigma\in\Aut\,\mathbf{S}^\mathrm{Pure} $ preserves the equivalence relations (\ref{qak}) and   commutes with the action of $\GL_2(\A_f),$ that is (\ref{autg}) holds and, for $q\in \GL_2(\Q),$ 
\be \label{aut0} q\tau=q\Longrightarrow [\tau, qa]^\sigma=[\tau,a]^\sigma.\ee 
Also, condition (iii) of \ref{definitionS} imposes the condition for arbitrary $\tau_1,\tau_2\in \HH$ and arbitrary $a_1,a_2\in \GL_2(\A_f),$
$$  [\tau_1,a_1], [\tau_2,a_2]\in \Ss^\mu \Rightarrow
  [\tau_1,a_1]^\sigma, [\tau_2,a_2]^\sigma\in \Ss^{\mu^\sigma}.$$
This implies for $E$ and $E_K$ defined in \ref{Fact2.6}(iii)
\be \label{components} E_K(s_1,s_2) \Leftrightarrow E_K(s_1^\sigma,s_2^\sigma);\ \ \ E(s_1,s_2) \Leftrightarrow E(s_1^\sigma,s_2^\sigma).\ee

On $\Ss\setminus \bigcup_{\tau\in \mathrm{CM}}\Ss_\tau$ the action of 
$\GL_2(\A_f)/Z(\Q)$ is free and hence the only condition on automorphisms of the $\GL_2(\A_f)/Z(\Q)$-set is the condition (\ref{components}).
 
 So we concentrate on the case that $\tau$ is quadratic, in which case
 $\Ss_\tau= \Ss_{\sqrt{-n}}$ for some 
 positive square-free integer $n$.

Note that   $\Ss_\tau$ is invariant under   $\sigma \in \Aut\, \mathbf{S}^\mathrm{Pure}.$ Indeed, $s\in \Ss_\tau$ if and only if $s=a\inv *\tau$ for some $a\in \GL_2(\A_f),$ which is equivalent by \ref{specpoints}(i) to the condition that $s$ is  fixed by a $g\in a\inv G_\tau a.$ But $g*s=s$  implies $g*s^\sigma=s^\sigma$  by (\ref{autg}). 

 Also, by \ref{specpoints} $\Ss_\tau=\Ss_\tau(\mathrm{CM}),$ that is the orbit consists of special points.
 
 We study $\mathbf{S}^\mathrm{Pure}_\tau,$ the substructure of $\mathbf{S}^\mathrm{Pure}$ on the set   $\Ss_\tau,$ $\tau$ quadratic.

\epk
\bpk \label{rtsProp} {\bf Proposition.}{\em For a quadratic $\tau,$ for any
$\sigma\in \Aut\, \mathbf{S}^\mathrm{Pure}_\tau$ 
there is 
$r=r_\tau(\sigma)\in N_\tau $ such that 
 \be\label{rsig} [\tau, a]^\sigma:= [\tau, r a],\mbox{ for all }a\in \GL_2(\A_f). \ee

Conversely, every $r\in  N_\tau $ determines a $\sigma=\sigma_{r,\tau}\in  \Aut\, \mathbf{S}^\mathrm{Pure}_\tau$ by the formula {\rm (\ref{rsig})} and 
\be \label{rsig2} \sigma_{r,\tau}=\sigma_{r',\tau} \Leftrightarrow
r'\cdot r\inv\in G_\tau.
\ee
}

 {\bf Proof.}
Suppose 
$$[\tau, 1]^\sigma=[\tau, r_\tau(\sigma)], \mbox{ for some } r_\tau(\sigma)\in \GL_2(\A_f),$$
and hence, for any $a\in \GL_2(\A_f),$
 \be\label{rts} [\tau, a]^\sigma=[\tau, 1\cdot a]^\sigma=a\inv\cdot [\tau,1]^\sigma= [\tau, r_\tau(\sigma)\cdot a].\ee
 
\medskip

We claim that  $r_\tau(\sigma)\in N_\tau .$
Indeed, let $\tau$ be the fixed point of $q\in  \GL_2^+(\Q),$ $q\notin Z(\Q).$ Then by \ref{specpoints} $[\tau,1]$ is a fixed point of $q$ and $q\inv$ and so 
$$[\tau,q]=[\tau,1]\mbox{ and } [\tau, r_\tau(\sigma)\cdot q]=  
[\tau,r_\tau(\sigma)].$$ 

Hence $[\tau,r_\tau(\sigma)]$ is a fixed point of $q$ and thus by \ref{specpoints}(ii) $r_\tau(\sigma)\in N_\tau .$

Conversely, if $r\in N_\tau$ then for each $q\in G_\tau\setminus Z(\Q)$ there is $q'\in G_\tau,$ 
$$[\tau, r]=[\tau, q'r]=[\tau,rq]=q* [\tau,r].$$ 
 
Define  the map $\sigma_{r,\tau}: \Ss_\tau\to \Ss_\tau,$  a partial transformation of $\Ss,$   by setting
for each $a\in \GL_2(\A_f)$ 
$$[\tau,a]^{\sigma_{r,\tau}}:= [\tau,ra].$$
Then $\sigma_{r,\tau}$ trivially satisfies (\ref{autg}). Also  (\ref{aut0}) is satisfied   because $rq=q'r$. This proves that $\sigma_{r,\tau}$ acts as an automorphism on $\Ss_\tau$ and $r=r_\tau(\sigma_{r,\tau}).$

Now note that by (\ref{qak})
$[\tau, r]=[\tau,1],$ for $r\in \GL_2(\A_f),$ holds  if and only if 
$r\in G_\tau.$ This proves (\ref{rsig2}). 
 $\Box$

 \epk

\bpk \label{12.8} 
    \label{AutC/F} {\bf Proposition} (The action of $\Aut\,\C$ on $\Ss_{\tau}$.)

 {\em  Let $\tau:=\sqrt{-n}$ and $\sigma\in \Aut(\C).$ }
 
(i) {\em $\sigma$ acts on $\Ss_{\tau}$ exactly by the automorphisms $$r=r_{\tau}(\sigma):\ [\tau,a]\ \mapsto\  [\tau, r a];\ \ r\in N_{\tau} /G_{\tau}.$$ }
  
  (ii)
$$N_{\tau} /G_{\tau}\cong\Aut\, \Ss_{\tau}^\mathrm{Full}=\Aut\, \Ss_{\tau}^\mathrm{Pure}\cong \Gal(\Q(\tau)^\mathrm{ab}:\Q)$$

 (iii) For any $g\in \GL_2(\A_f),$
 $$g^\sigma=g.$$
 
 (iv) $$\Aut \mathbf{S}^\mathrm{Full}\subseteq \Aut \mathbf{S}^\mathrm{Pure}.$$

{\bf Proof.}
Compare the conclusion of \ref{rtsProp} with \cite{Milne0}, 12.8 describing the action on the special points $\Ss_\tau$  of the canonical model by
$\sigma\in \Gal(\Q(\tau)^\mathrm{ab}:\Q(\tau)).$ 

In the notation of \cite{Milne0}, (60)-(62):

$\tau=x,$ $\Q(\tau)=E(x)$ and 
$T(\A_f)=T_\tau $ as in our \ref{S+}. Further
\be \label{GT1} r_x: \A^\times_{\Q(\tau)}\ \ \to \  T(\A_f); \ \ a+b\tau\ \ \mapsto \ \ \left( \begin{array}{ll}
\ a\ nb\\ -b\ \ a
\end{array}\right).\ee
 
 Consequently, $r_\tau(\sigma)$ is the same as $r_x(s)$   of \cite{Milne0}, 12.8, for $s\in \A^\times_{\Q(\tau)}$ such that $\sigma=\mathrm{art}(s),$ the Artin symbol. 
 Thus the matrices $r_x(s)\in T_\tau $  of (\ref{GT1}) deliver all the Galois automorphism of $\Ss_\tau(\mathrm{CM})$ which fix
 $\Q(\tau).$ 
 
 Now let $\mathsf{i}: \C \to \C$ be  complex conjugation.
 $\mathsf{i}$ induces a non-trivial automorphism of $\Q(\tau)$ and an automorphism of $\Ss_\tau(\mathrm{CM}).$ 
 
Claim.  $\mathsf{i}$
  corresponds to the action by an element
$$r_{\tau}(\mathsf{i})= \left( \begin{array}{ll}
-c\ nd\\ -d\ \ c
\end{array}\right)\in N_{\tau}\setminus T_{\tau}.$$
  
 Proof of Claim. The transformation $\eta: \HH\cup \HH^-\to \HH\cup \HH^-$
 induced on $\HH\cup \HH^-\subset \Ss$ by $\mathbf{i}$ is studied (in much more general setting) in
 \cite{MilneCC}. Lemma 3.2  of that paper states in particular
 that (in our notation $n:=\rrr_{-1}$), for $g\in \GL_2(\Q),$    
 $$\eta([\tau,g])=[\tau,\rrr_{-1}\cdot g].$$  
In particular, $$\eta(\tau)=\rrr_{-1}\cdot\tau.$$

It implies, as proved further in \cite{MilneCC}, Conjecture B (a form of Langlands' conjecture on complex conjugation): 
$$[\tau,g]^\mathbf{i}=[\eta(\tau),g]=[\tau,\rrr_{-1}\cdot g]$$  
for all $g\in \GL_2(\A_f).$
 
 \medskip
 
Now let
  $\rho\in \Gal(\Q(\tau)^\mathrm{ab}:\Q).$
  Then, either it fixes  $\Q(\tau)$ and has the form (\ref{GT1}) or $\mathsf{i}\circ\rho$  fixes  $\Q(\tau)$ and so $\rho$ can be represented by the product of a matrix from 
  $T_{\tau}$ with the above  matrix from   $N_{\tau},$ which gives us again a matrix of the form $\left( \begin{array}{ll}
-c\ nd\\ -d\ \ c
\end{array}\right)$ from  $N_{\tau}.$ 
This concludes the proof of (i) and (ii).

Proof of (iii). Since
$$ (g*[\tau,a])^\sigma=([\tau,ag\inv])^\sigma=[\tau, rag\inv]=g*[\tau,a]^\sigma,$$
the set of special points in $C_g$ (the graph of $g$ on $\Ss$) is invariant under $\sigma.$ 
Hence the set of special points in the algebraic curve $C_{g,N}\subset \Y(N)\times \Y(N)$  is invariant under $\sigma.$ But the set of special points is Zariski dense in the special curve $C_{g,N},$ and hence $C_{g,N}$ is invariant under $\sigma.$ Finally, this implies $C_g,$ the limit of the $C_{g,N},$ is invariant under $\sigma.$

Proof of (iv). We have established in (iii) that $\mathbf{S}^\mathrm{Full}$ satisfies (\ref{autg}), that is $s\mapsto g*s$ is definable in the structure. We conclude  by \ref{Fact2.6}
that
$\mathbf{S}^\mathrm{Full}$ is an expansion of $\mathbf{S}^\mathrm{Pure}$ by possibly extra relations, which implies the embedding of the automorphism groups.
$\Box$

    \epk

 \bpk \label{actHn} {\bf The action of 
 $\Aut\, \C$ on canonical models over $\HH^n,$
 the direct product of the upper half-plane and the action of $\Aut\ \mathbf{S}^\mathrm{Pure}(\mathrm{CM})$ on $\mathbf{S}^\mathrm{Pure}(\mathrm{CM})$}

Recall the definitions  (see \ref{introSh}) related to the canonical models of the  Shimura variety corresponding to $( \GL_2^n,(\HH\cup -\HH)^n )$:
$$\Ss^n_K:= \GL_2(\Q)^n\backslash  (\HH\cup -\HH)^n \times \GL_2^n(\A_f)/K,$$
where $K$   is a compact open subgroup of  $\GL_2^n(\A_f).$

The elements of  $\Ss^n$ can be written as
  $$[\bar{\tau}, a],\mbox{ for } \bar{\tau}\in \HH^n,\ a\in  \GL_2( \A_f)^n$$  
and the condition that co-ordinates of $\bar{\tau}$ are quadratic points implies that $\bar{\tau}$ is a fixed point of a $\bar{q}\in  \GL_2(\Q)^n.$ 

Consequently,  a $\sigma\in \Aut\,\C$ acts on  $[\bar{\tau}, a]$
so that
\be\label{rbart} [\bar{\tau}, a]^\sigma=[\bar{\tau}, \bar{r}_\tau(\sigma) a]\ee
where \be \label{rbart2}\bar{r}_{\bar{\tau}}(\sigma)=\la r_{\tau_1}(\sigma),\ldots,r_{\tau_n}(\sigma)\ra\mbox{ and }r_{\tau_i}(\sigma)\in N_{\tau_i},\ i=1,\ldots,n.\ee
When $\sigma\in \Aut(\C/\Q(\bar{\tau}))$ then  $$\bar{r}_{\bar{\tau}}(\sigma)\in T_{\bar{\tau}} ,$$
 where $T_{\bar{\tau}}< \GL_2^n$ is the  torus over $\Q$ for the special point $\bar{\tau},$ see \cite{Milne0}, Ch.12, (60) and (61). More generally, $\bar{r}_{\bar{\tau}}(\sigma)\in N_{\bar{\tau}} $ for  
 $\sigma\in \Aut\,\C.$

\medskip

When $\bar{\tau}$ is a special point, then 
(\ref{rbart}) and (\ref{rbart2}) are equally applicable to  $\sigma\in \Aut\ \mathbf{S}^\mathrm{Pure}(\mathrm{CM})$ due to \ref{AutC/F}.

\medskip

{\bf Remark.} An important property of $\sigma$ is the {\bf preservation of connected components}, see \ref{Fact2.6}.
 Hence, for some $\mu=\mu(\sigma)\in \hat{\Z}^\times$
\be \label{propE}\bigwedge_{i=1}^n\, [\tau_i, r_{\tau_i}(\sigma)]\in \Ss^\mu
\ee 
\epk
\bpk
Let
$$\Aut\, \mathbf{S}^\mathrm{Pure}_{\bar{\tau}}$$
be the automorphism group of the pure structure on $\Ss_{\bar{\tau}}=\bigcup_{i=1}^n \Ss_{\tau_i}.$ By the remark \ref{2.5+} following the definition of $\mathbf{S}^\mathrm{Pure},$ $\sigma$ acts on $\Ss_{\bar{\tau}}$ preserving the action of $\GL_2(\A_f)$ and the equivalence $E.$ The first is the condition (\ref{autg}) and the second (\ref{propE}). Since we represent the points of $\Ss$ as double-cosets we also need (\ref{aut0}).
  
That is a transformation
$\sigma: \Ss_{\bar{\tau}}\to \Ss_{\bar{\tau}}$ is in $\Aut \mathbf{S}^\mathrm{Pure}_{\bar{\tau}}$ if and only if $\sigma$ satisfies (\ref{autg}),(\ref{aut0}) and (\ref{propE}) for points in $\Ss_{\bar{\tau}}.$

We are going to study the case where $\bar{\tau}$ is a special point of 
$\HH^n,$ that is the $\tau_i$ are quadratic imaginary. In this case
$$\Ss_{\bar{\tau}}=\Ss_{\bar{\tau}}(\mathrm{CM}).$$

We denote
$$\Aut\, \mathbf{S}^\mathrm{Full}_{\bar{\tau}},$$
the automorphism group of the full structure on $\Ss_{\bar{\tau}}.$ 
Since the structure on CM-points of the full structure is determined by the embedding into $\Q^\mathrm{alg},$ the algebraic points of the field, 
 $$\Aut\, \mathbf{S}^\mathrm{Full}_{\bar{\tau}}=\{ \sigma\in \Aut\, \mathbf{S}^\mathrm{Pure}_{\bar{\tau}}| \ \sigma \mbox{ induced by a Galois automorphism over }\Q\}$$
 and so $\Aut\, \mathbf{S}^\mathrm{Full}_{\bar{\tau}}$ can be naturally identified with a Galois group over $\Q.$
\epk

\bpk \label{314}
We will work with the following groups which are sub-direct products  of  groups classfied in \ref{AutC/F}: 

For $\bar{\tau}:=\la \tau_1,\ldots,\tau_n\ra= \la\sqrt{-m_1},\ldots,\sqrt{-m_n}\ra,$
$$T^*_{\bar{\tau}}:= \left\lbrace \la r_{1},\ldots,r_{n}\ra\in
T_{\tau_1}\times \ldots \times T_{\tau_n}: \det r_{1}=\ldots = \det r_{n}\right\rbrace$$
$$T^0_{\bar{\tau}}:= T^0_{\tau_1}\times \ldots \times T^0_{\tau_n},\ \ \ \ T^0_{\tau_i}:=\{ r_i\in T_{\tau_i}: \ \det r_i=1\} $$
$$G_{\bar{\tau}}:=  G_{\tau_1}\times \ldots \times G_{\tau_n}$$
and
$$\bar{T}^*_{\bar{\tau}}:=T^*_{\bar{\tau}}\cdot G_{\bar{\tau}}/ G_{\bar{\tau}}, \ \
\bar{T}^0_{\bar{\tau}}:=T^0_{\bar{\tau}}\cdot G_{\bar{\tau}}/G_{\bar{\tau}}.$$

\medskip

We say, for $m_1,\ldots,m_n$ square-free positive distinct integers that 
$\sqrt{-m_1},\ldots,\sqrt{-m_n}$ are {\bf  independent over $\Q$} if 
 $|\Q(\sqrt{-m_1},\ldots,\sqrt{-m_n}):\Q|=2^n,$
  equivalently, $$\sqrt{-m_{i+1}}\notin  \Q(\sqrt{-m_1},\ldots,\sqrt{-m_i}),\mbox{ for }i=1,\ldots,n-1.$$
  

\medskip

For  $\sigma\in \Aut\, \mathbf{S}^\mathrm{Pure}_{\bar{\tau}},$ the notation  $\bar{r}_{\bar{\tau}}(\sigma)$ is defined by (\ref{rbart}) and (\ref{rbart2}).
\epk
Our {\bf assumption} throughout the rest of the paper is that 
  $$\bar{\tau}=\la \tau_1,\ldots,\tau_n\ra=\la \sqrt{-m_1},\ldots,\sqrt{-m_n}\ra,$$    independent over $\Q.$ 
\bpk \label{TT1} {\bf Proposition.}

(i) {\em  For all $\sigma\in \Aut\, \mathbf{S}^\mathrm{Pure}_{\bar{\tau}}:$

  $\bar{r}_{\bar{\tau}}(\sigma)\in \prod_i N_{\tau_i} $ and
   $\bar{r}_{\bar{\tau}}(\sigma)$ can be represented as
    $\la r_1,\ldots, r_n\ra=\bar{r}_{\bar{\tau}}(\sigma)$    
   so that, for some $q_1,\ldots, q_n\in \Q_+,$ 
$$q_1\cdot\det r_1=\ldots =q_n\cdot\det r_n=\mu(\sigma)\in \hat{\Z}^{\times}.$$

Moreover, $\mu(\sigma)=1$ if and only if all the components $\Ss^\lambda$ of $\Ss$ are invariant under $\sigma,$ and $\mu(\sigma)$
 runs through all values $\mu\in \hat{\Z}^{\times}$ as $\sigma$ runs in $\Aut\, \mathbf{S}^\mathrm{Full}_{\bar{\tau}}.$
}

\medskip

(ii) {\em For any $\sigma\in \Aut\, \mathbf{S}^\mathrm{Pure}_{\bar{\tau}}$
one can choose a representative $\la r_1,\ldots,r_n\ra$ of $\bar{r}_{\bar{\tau}}(\sigma)$ and $q\in \Q_+$
so that, for some partition $\{ 1,\ldots,n\}=\{ i_1,\ldots, i_k\}\dot{\cup} \{ i_{k+1},\ldots, i_n\},$
the following hold:
\be\label{riri2}  \begin{array}{lll}\det r_{i_1}=\ldots=\det r_{i_k} \mbox{ and }r_{i_1}\in T_{\tau_{i_1}},\ldots,r_{i_k}\in T_{\tau_{i_k}},\\
\det r_{i_{k+1}}=\ldots=\det r_{i_n} \mbox{ and }r_{i_{k+1}}\notin T_{\tau_{i_{k+1}}},\ldots,r_{i_n}\notin T_{\tau_{i_n}},\\
\det r_{i_1}=\ldots=\det r_{i_k}=q\det r_{i_{k+1}}=\ldots=q\det r_{i_n}
\end{array}
\ee
 }

(iii) 
{\em Given $\la r_1,\ldots,r_n\ra\in \prod_i N_{\tau_i} $ and   $q\in \Q_+$ suppose there exists  a partition $\{ 1,\ldots,n\}=\{ i_1,\ldots, i_k\}\dot{\cup} \{ i_{k+1},\ldots, i_n\}$ satisfying
(\ref{riri2}). 
Then 
  there exists $\sigma\in \Aut\, \mathbf{S}^\mathrm{Pure}_{\bar{\tau}}$
such that $\la r_1,\ldots,r_n\ra$ is a representative of $\bar{r}(\sigma).$ }

\medskip

{\bf Proof.} (i) By the preservation of connected components property (\ref{propE}), given $\sigma$ there is a $\mu$ such that 
for each $\tau_i:= \sqrt{-m_i},$
$$[\tau_i, r_{\tau_i}(\sigma)]\in \Ss^\mu, \mbox{ and hence } \det r_{\tau_i}(\sigma)=q_i\inv\cdot\mu, \ \mu\in \hat{\Z}^\times, \mbox{ some } q_i\in \Q_+.$$

This proves the first statement of (i). For the second statement
 use the fact that $\Gal(\Q^\mathrm{ab}:\Q)$ acts faithfully and transitively on the set of components $\Ss^\mu$ of $\Ss$ (see \cite{Milne1}, 3.4) so for each $\mu$ there is a $\sigma$ and $\bar{r}_{\bar{\tau}}(\sigma)$ such that  $\bigwedge_i \det r_{\tau_i}(\sigma)= q_i \mu.$

(ii) In matrix terms 
we have $$r_{\tau_i}(\sigma)=r_i=\left( \begin{array}{ll}\ a_i\ \ \ b_im_i\\ 
-\epsilon_i b_i\ \ \epsilon a_i
\end{array}\right),\ \ \ \det r_i=\epsilon_i  (a_i^2 + b_i^2m_i)=  q_i\inv\mu,\ ( a_i,b_i\in \A_f)$$
where $\epsilon_i$ is $+1$ if  $r_i\in T_{\sqrt{-m_i}}$  and $-1,$ if 
$r_i\in N_{\sqrt{-m_i}}\setminus T_{\sqrt{-m_i}}.$

As a result we have the equalities
\be\label{q1} \epsilon_i q_i (a_i^2+b_i^2m_i)= \mu = \epsilon_j q_j (a_j^2+b_j^2m_j);
\ i,j=1,\ldots,n.\ee

Assume that $\epsilon_i=\epsilon_j.$ 
Then the equality between the algebraic expressions takes place for 
some $a_i,b_i,a_j,b_j\in \R,$ not both $a_i,b_i$ zero, as well. By the Hasse principle, there are $s_i,t_i\in \Q,\ \ i=1,\ldots, n$
such that not both $s_i,t_i$ zero and
\be \label{q2}\epsilon_i q_i (s_i^2+t_i^2m_i)=\epsilon_j q_j (s_j^2+t_j^2m_j).\ee

Consider, for $i=1,\ldots,n,$ the  matrices with rational entries,
 $$g_i=\left( \begin{array}{ll}\ s_i\ t_im_i\\ 
 -t_i\ \  s_i
\end{array}\right),$$
which by definition belongs to $G_{\sqrt{-m_i}}.$

By dividing the algebraic expressions of (\ref{q1}) on respective expressions of (\ref{q2}) we get
$$(a_i^2+b_i^2m_i)(s_i^2+t_i^2m_i)\inv=(a_j^2+b_j^2m_j)(s_j^2+t_j^2m_j),$$
and recalling the formula for $\det r_i$  
we have  $$\det r_i g_i\inv=\epsilon_i(a_i^2+b_i^2m_i)(s_i^2+t_i^2m_i)\inv=\epsilon_j(a_j^2+b_j^2m_j)(s_j^2+t_j^2m_j)\inv=
\det r_jg_j\inv.$$
Choose  representatives for $\bar{r}_{\bar{\tau}}(\sigma)$ modulo $G_{\bar{\tau}}$ as 
$$r_i:= r_i g_i\inv$$

Suppose $\epsilon_{i_1}=\ldots=\epsilon_{i_k}=1$ and $\epsilon_{i_{k+1}}=\ldots=\epsilon_{i_n}=-1.$

We get, for some $q\in \Q_+,$ from above that
$$\begin{array}{lll}\det r_{i_1}=\ldots=\det r_{i_k} \mbox{ and }r_{i_1}\in T_{i_1},\ldots,r_{i_k}\in T_{i_k}\\
\det r_{i_{k+1}}=\ldots=\det r_{i_n} \mbox{ and }r_{i_{k+1}}\notin T_{i_{k+1}},\ldots,r_{i_n}\notin T_{i_n}\\
\det r_{i_1}=\ldots=\det r_{i_k}=q\det r_{i_{k+1}}=\ldots=q\det r_{i_n}
\end{array}
$$ 

(iii) The condition $\bar{r}=\la r_1,\ldots,r_n\ra\in \prod_i N_{\tau_i} $
implies that $\bar{r}$ determines a transformation of $\Ss_{\bar{\tau}}$ which satisfies (\ref{autg}) and (\ref{aut0}). The conditions (\ref{riri2}) implies that $\bar{r}$ satisfies (\ref{propE}). Thus the transformation belongs to $\Aut\, \mathbf{S}^\mathrm{Pure}_{\bar{\tau}}.$
$\Box$
\epk 
\bpk\label{Cor316} {\bf Corollary.} {\em Let $$N^*_{\bar{\tau}} =\{ \la r_1,\ldots,r_n\ra\in \prod_{i=1}^n N_{\tau_i} : \ \exists_{i=1}^n 
q_i\in \Q_+, \exists \mu\in \hat{\Z}^{\times}\\ \bigwedge_{i=1}^n 
q_i\cdot\det r_i=\mu\}$$
and
$$\bar{N}^*_{\bar{\tau}} := N^*_{\bar{\tau}} /G_{\bar{\tau}}.$$
Then
$$\bar{r}_{\bar{\tau}}/G_{\bar{\tau}}: \sigma\mapsto \bar{r}_{\bar{\tau}}(\sigma)/G_{\bar{\tau}}$$
is an isomorphism
$$\Aut\, \mathbf{S}^\mathrm{Pure}_{\bar{\tau}}\cong \bar{N}^*_{\bar{\tau}} $$
as groups of transformations of $\Ss_{\bar{\tau}}.$
 }
\epk

\bpk Set 
$$\bar{N}^{\Gal}_{\bar{\tau}} $$
to be the image of the Galois group in $\bar{N}^*_{\bar{\tau}} $
under the canonical map
$$\sigma\mapsto \bar{r}_{\bar{\tau}}(\sigma)/G_{\bar{\tau}}.$$

Equivalently,
$$N^{\Gal}_{\bar{\tau}} =\{ \bar{r}\in N^*_{\bar{\tau}} : \ \exists \sigma\in \Gal(\Q(\mathrm{CM}):\Q), \ 
\bar{r}/G_{\bar{\tau}}=\bar{r}_{\bar{\tau}}(\sigma)/G_{\bar{\tau}}\}$$
and $$\bar{N}^{\Gal}_{\bar{\tau}} := N^{\Gal}_{\bar{\tau}} /G_{\bar{\tau}}.$$

\medskip

Let $$\mathrm{CM}_{\tau_i}=\bigcup_{N=1}^\infty j_N(\Ss_{\tau_i})$$ be  the CM-points on all the curves $\Y(N)$ corresponding to the imaginary quadratic $\tau_i.$
 Set $$\mathrm{CM}_{\bar{\tau}}:=\bigcup_{i=1}^n\mathrm{CM}_{\tau_i}=\bigcup_{N=1}^\infty j_N(\Ss_{\bar{\tau}}).$$
\epk
 We denote $$\F_{\bar{\tau}}=\Q(\tau_1)^\mathrm{ab}\cdot \ldots \cdot \Q(\tau_n)^\mathrm{ab},$$
the composite of the $n$ fields.

\bpk \label{318-}\label{318} {\bf Lemma.} 

{\em
$$\Q(\mathrm{CM}_{\bar{\tau}})=\F_{\bar{\tau}}.$$
and
 $$\Gal(\F_{\bar{\tau}}: \Q)\cong \Aut\, \mathbf{S}^\mathrm{Full}_{\bar{\tau}}\cong \bar{N}^{\Gal}_{\bar{\tau}} $$
as groups of transformations of $\Ss_{\bar{\tau}}.$}

{\bf Proof.}   By \ref{AutC/F}(ii)
 $\Q(\mathrm{CM}_{\tau_i})=\Q(\tau_i)^\mathrm{ab}.$ Hence
 $\Q(\bigcup_{i=1}^n\mathrm{CM}_{\tau_i})=\F_{\bar{\tau}}.$

Elements $\sigma$ of $\Aut\, \mathbf{S}^\mathrm{Full}_{\bar{\tau}}$ by construction act on  $\bigcup_i\mathrm{CM}_{\tau_i}$   as an automorphism group of the field, hence 
$$\Aut\, \mathbf{S}^\mathrm{Full}_{\bar{\tau}}\cong \Gal(\F_{\bar{\tau}}: \Q).$$ 

By \ref{Cor316} and definitions  $\bar{N}^{\Gal}_{\bar{\tau}} $ is isomorphic to the subgroup of the group of transformations 
 $\Aut\, \mathbf{S}^\mathrm{Pure}_{\bar{\tau}}$ induced by Galois automorphisms of $\Q(\mathrm{CM}_{\bar{\tau}})$ over $\Q,$ which is exactly $\Aut\, \mathbf{S}^\mathrm{Full}_{\bar{\tau}}.$ This proves the second isomorphism.
$\Box$
\epk
 Let   $$\bar{T}^{\Gal}_{\bar{\tau}}=\bar{T}_{\bar{\tau}}\cap \bar{N}^{\Gal}_{\bar{\tau}} $$
 and
 $$\bar{T}^{0,\Gal}_{\bar{\tau}}=\bar{T}^0_{\bar{\tau}}\cap \bar{N}^{\Gal}_{\bar{\tau}} .$$
\bpk \label{(b)} {\bf Lemma.} {\em $$\bar{T}^{0,\Gal}_{\bar{\tau}}\cong \Gal(\F_{\bar{\tau}}:\Q^\mathrm{ab})$$
and $\bar{T}^{0,\Gal}_{\bar{\tau}}$ is isomorphic to an open  subgroup of $\bar{T}^{0}_{\bar{\tau}}:$

$$\bar{T}^{0,\Gal}_{\bar{\tau}}\le_\mathrm{open}\bar{T}^0_{\bar{\tau}}.$$}

{\bf Proof.}  First we claim that as groups of transformations, for each $i,$
$$\Gal(\Q(\tau_i)^\mathrm{ab}:\Q^\mathrm{ab})\cong \bar{T}^{0,\Gal}_{\tau_i} $$ via
 the canonical map
$$\sigma\mapsto \bar{r}_{\bar{\tau}}(\sigma).$$
In order to prove it note that  
$T^0_{\tau_i}\cong \A^0_{ \Q(\tau_i)}$ where  $\A^0_{ \Q(\tau_i)}$
is the  group of finite  $\Q(\tau_i)$-adeles of $\Q$-norm 1.
Now    
the isomorphism $$\A^0_{ \Q(\tau_i)}/\pm 1\cong \Gal(\Q(\tau_i)^\mathrm{ab}:\Q^\mathrm{ab}),$$
a corollary to Artin's reciprocity law, completes the proof of our claim.

Now the definition of $\bar{T}^{0,\Gal}_{\bar{\tau}} $ 
together with the second statement of \ref{318-} and the claim above gives us
$$\bar{T}^{0,\Gal}_{\bar{\tau}}\cong  \left( \Gal(\Q(\tau_1)^\mathrm{ab}:\Q^\mathrm{ab})\times \ldots \times \Gal(\Q(\tau_n)^\mathrm{ab}:\Q^\mathrm{ab})\right)\cap \Gal(\F_{\bar{\tau}}:\Q).$$
It is easy to check using the definition of  $\F_{\bar{\tau}}$
that the right-hand side of the latter is exactly $\Gal(\F_{\bar{\tau}}:\Q^\mathrm{ab}).$

Thus
$$\bar{T}^{0,\Gal}_{\bar{\tau}}\cong \Gal(\F_{\bar{\tau}}:\Q^\mathrm{ab}).$$

For the proof of the second statement of the Lemma we use the proof of the adelic Mumford-Tate conjecture for the product of CM elliptic curves, \cite{Cadoret--Moonen} and \cite{Ullmo--Yafaev}.

Let ${\bf T}^*_{\bar{\tau}}$ denote the $\mathbb{Q}$-torus such that ${\bf T}^*_{\bar{\tau}}(\mathbb{A}_f)=T^*_{\bar{\tau}}$. We can form a CM-Shimura datum $({\bf T}*_{\bar{\tau}},\{h\})$ associated with a product $E_1\times\cdots\times E_n$ of elliptic curves, with $E_i$ having CM by $\mathbb{Q}(\tau_i)$. As in \cite{Ullmo--Yafaev},  Remark 2.8, for any neat compact open subgroup $K$ of $T^*_{\bar{\tau}}$ and $L$ a finite extension of $\mathbb{Q}(\bar{\tau})$ such that $[h,1]\in{\rm Sh}_{K}({\bf T}^*_{\bar{\tau}},\{h\})(L)$, we obtain a homomorphism
\[\rho:{\rm Gal}(\bar{L}/L)\to K,\]
which describes the action of ${\rm Gal}(\bar{L}/L)$ on the fiber of
\[{\rm Sh}({\bf T}^*_{\bar{\tau}},\{h\})\to{\rm Sh}_{K}({\bf T}^*_{\bar{\tau}},\{h\})\]
over $[h,1]$ (which is a $K$-torsor). In particular, $\rho$ factors through the quotient ${\rm Gal}(\F_{\bar{\tau}}:L)$.

As explained on the following page of \cite{Ullmo--Yafaev}, $\rho$ yields the ${\rm Gal}(\bar{L}:L)$ representation on the Tate-module of $E_1\times\cdots\times E_n$. Take $L$ to satisfy $ \Q(\bar{\tau})\subseteq L\subset\F_{\bar{\tau}}.$
By \cite{Cadoret--Moonen}, Theorem A (ii), provided $h_\mathbb{R}$ is ``maximal'' (see Definition 2.1), the latter has open image in $T^*_{\bar{\tau}}$ and, therefore,
\[{\rm Gal}(\F_{\bar{\tau}}:L)\cap T^0_{\bar{\tau}}={\rm Gal}(\F_{\bar{\tau}}:L\mathbb{Q}^{\rm ab})\]
(the equality due to the fact that $T^0_{\bar{\tau}}$ is the kernel of the determinant map $T^*_{\bar{\tau}}\to\mathbb{A}^\times_f$) is open in $T^0_{\bar{\tau}}$.

To see that $h_\mathbb{R}$ is maximal, let $\sigma_i$ and $\bar{\sigma}_i$ denote distinct embeddings of $\mathbb{Q}(\tau_i)$. The cocharacter group of ${\bf T}^*_{\bar{\tau}}$ is equal to the set of elements 
\[\sum_{i=1}^n a_i\sigma_i+b_i\bar{\sigma}_i\]
such that $a_i+b_i$ is independent of $i$. It has a basis consisting of the $\sigma_i-\bar{\sigma}_i$, for $i=1,\ldots,n$, and $\mu$, where $\mu=\sum_i\sigma_i$ is the cocharacter associated with $h$. The ${\rm Gal}(\bar{\mathbb{Q}}:\mathbb{Q})$-orbit of $\mu$ consists of the sums of $n$ elements in which, for each $i$, the sum includes either $\sigma_i$ or $\bar{\sigma}_i$. Since 
\[\sigma_i-\bar{\sigma}_i=\mu-(\Sigma_{j\neq i}\sigma_j+\bar{\sigma}_i),\] the result follows from the the paragraph following Definition 2.3 of
\cite{Cadoret--Moonen}. $\Box$

\epk

\bpk \label{320} {\bf Lemma.} 
{\em The following diagram commutes 
$$\xymatrix{1 \ar[r] &\bar{T}^{0,\Gal}_{\bar{\tau}}  \ar[r]^{i}  \ar[d]^{i}&\bar{N}^{\Gal}_{\bar{\tau}}   \ar[r]^{\ \ \ \det/\Q_+}\ar[d]^{i} &\hat{\Z}^\times \ar[d]^{\mathrm{id}}\ar[r] &1 \\ 1 \ar[r] &\bar{T}^0_{\bar{\tau}}\ \ar[r]^{i} & \bar{N}^*_{\bar{\tau}} \ar[r]^{\ \det/\Q_+} &\hat{\Z}^\times \ar[r] &1}$$
and both lines are exact. The arrows marked $i$ are natural embeddings, and 
$\det/\Q_+$ applied to an element represented by $\la r_1,\ldots,r_n\ra$
is $\det r_i$ modulo $\Q_+,$ which is of the same value for all $i.$}

{\bf Proof.} $\det/\Q_+(r_1,\ldots,r_n)$ takes its values in $\A_f^\times/\Q_+^\times\cong \hat{\Z}^\times.$  
It
 is well-defined by the definitions of 
$\bar{N}^{\Gal}_{\bar{\tau}} $ and $\bar{N}^*_{\bar{\tau}} ,$
and is surjective on both lines by \ref{TT1}(i). 
The kernel of $\det/\Q_+$ consists of tuples $\bar{r}$ of matrices $r_i=\left( \begin{array}{ll}
a_i\ \ b_im_i\\
-b_i\ a_i
\end{array}\right)$ such that $\det r_i=a_i^2+m_ib_i^2=q_i\in \Q^+,$ where $a_i,b_i\in \A_f.$ By the Hasse principle there exists $s_i,t_i\in \Q$ such that $s_i^2+m_it_i^2=q_i$ and thus there exists $g_i\in G_{\tau_i}$ such that $\det r_ig_i\inv=1.$ Thus $\bar{r}\in   \bar{T}^0_{\bar{\tau}}.$

It follows that the kernel
is $\bar{T}^{0,\Gal}_{\bar{\tau}}$ in the top line and
$\bar{T}^0_{\bar{\tau}}$ in the bottom line. $\Box$  

\epk
\bpk \label{Cor321} {\bf Corollary.} {\em
$$\bar{N}^{\Gal}_{\bar{\tau}} \le_\mathrm{open} \bar{N}^{*}_{\bar{\tau}} $$
and
$$\Gal(\F_{\bar{\tau}}: \Q)\cong\Aut\, \mathbf{S}^\mathrm{Full}_{\bar{\tau}}\le_\mathrm{open} \Aut\, \mathbf{S}^\mathrm{Pure}_{\bar{\tau}}$$
as the transformation groups on $\Ss_{\bar{\tau}}.$}

{\bf Proof.} The first statement is immediate from \ref{(b)} -  \ref{320}. The second statement follows from the first together with \ref{Cor316} and \ref{318}. $\Box$
\epk
\bpk \label{N/T} {\bf Lemma.} {\em 
$$\bar{N}^\Gal_{\bar{\tau}} /\bar{T}^\Gal_{\bar{\tau}}  \cong \mathrm{C}_2^n,$$
where $\mathrm{C}_2$ is the 2-element group generated by the involution $\rrr_{-1}.$
The group acts on $\bar{T}^*_{\bar{\tau}} $ coordinate-wise by conjugation by  $\rrr_{\pm 1}$ on the $i$-th coordinate.
}

{\bf Proof.} By  \ref{318} and \ref{Cor321} $\bar{N}^\Gal_{\bar{\tau}} $ 
acts on the compositum $\F_{\bar{\tau}}$ of the fields $\Q(\tau_i)^\mathrm{ab}$ as its
Galois group over $\Q.$  At the same time coordinate-wise  $\bar{r}=\la r_1,\ldots,r_n\ra\in N^*_{\bar{\tau}}$ acts on
the $i$-th coordinate by $[\tau_i,a]\mapsto [\tau_i, r_i a],$ which by  \ref{12.8} corresponds to the Galois action on $\Q(\tau_i)^\mathrm{ab}.$ Moreover, $r_i(\sigma) \in T_{\tau_i} $ if and only if $\sigma$ fixes $\tau_i$ in $\Q(\tau_i)^\mathrm{ab}.$  If $r_i(\sigma) \notin T_{\tau_i} $ then $r_i(\sigma)$ acts on $T_{\tau_i} $ by conjugation as $\rrr_{-1}.$

Let $\delta\in \{ 1,-1\}^n$ be an arbitrary sequence with values $\pm 1.$
The condition of independence of  the $\tau_i$ guaranties the existence of $\sigma_\delta\in \Gal(\F_{\bar{\tau}}:\Q)$ such that 
$$\sigma_\delta(\tau_i)=\tau_i \Leftrightarrow \delta(i)=1.$$
Thus
$$ r_i(\sigma_\delta)\in T_{\tau_i} \Leftrightarrow \delta(i)=1.$$
It implies that 
 $$\bar{r}(\sigma_\delta)/T^\Gal_{\bar{\tau}} =\bar{r}(\sigma_{\delta'})/T^\Gal_{\bar{\tau}} \Leftrightarrow \delta=\delta'.$$
 $\Box$
\epk

\bpk \label{TT3}  
Now we summarise the facts proven above into a theorem characterising the Galois action on all the special points.

\medskip

Let $\F$ be
the composite of all $\Q(\tau)^\mathrm{ab}$ for all imaginary quadratic  $\tau.$

Let $P\subset \N$ be the set of primes, including 1. Clearly,
$\sqrt{P}:=\{ \sqrt{-p}:\ p\in P\} $ is a set of imaginary quadratics independent 
over $\Q$ (see \ref{314}) and generate the field containing all the quadratic imaginary extensions of $\Q.$

In the definitions below we take limits over finite tuples $\bar{\tau}$ in $\sqrt{P} .$
\medskip



$$\bar{T}^\Gal_* =\varprojlim \bar{T}^\Gal_{\bar{\tau}} ,$$
$$\bar{T}^{0,\Gal}_* =\varprojlim \bar{T}^{0,\Gal}_{\bar{\tau}} ,$$
$$\bar{N}^\Gal_* =\varprojlim \bar{N}^\Gal_{\bar{\tau}} .$$
$$\bar{N}^*_* =\varprojlim \bar{N}^*_{\bar{\tau}} .$$


Below we will claim that $\bar{N}^\Gal_* $ is {\bf a an "almost equal" subgroup  of } $\bar{N}^*_* ,$ and write $$\bar{N}^\Gal_*   \le_\approx\bar{N}^*_* ,$$
 meaning that $\bar{N}^\Gal_{\bar{\tau}} $ is a finite index subgroup of $\bar{N}^*_{\bar{\tau}} ,$ for each finite tuple $\bar{\tau}$ and both are projective limits along the same projective system indexed by finite tuples of elements of $\sqrt{P}.$ 

\epk
\bpk\label{MainGal}\label{GalCM}  {\bf Theorem.}

 
(i) {\em
 $$\bar{N}^\Gal_* \cong \Gal(\F: \Q)= \Gal(\Q(\mathrm{CM}): \Q)$$
and 
$$\bar{N}^\Gal_*   \le_\approx\bar{N}^*_* $$
as groups of transformations. }

(ii) {\em There is a short exact sequence of groups 
 $$1\to \bar{T}^\Gal_* \to \bar{N}^\Gal_*  \to \mathrm{C}_2^\omega \to 1$$
 where 
 $\mathrm{C}_2^\omega$ is the infinite product $\prod_{p\in P} \mathrm{C}_2(p)$ of 2-element groups isomorphic to
$\la \rrr_{-1}\ra$  acting on the infinite product of $2\times 2$ matrix group coordinate-wise by
 $$\rrr_{-1}: \left( \begin{array}{ll}
 \ a\ pc\\
 -c\ \ a
 \end{array}\right) \mapsto \left( \begin{array}{ll}
  a\ -pc\\
 c\ \ \ \ \ a
 \end{array}\right).$$}
 
 (iii) {\em There is a short exact sequence of groups $$1\to \bar{T}^{0,\Gal}_* \to \bar{N}^\Gal_*  \to \hat{\Z}^\times \to 1,$$
 $$\bar{T}^{0,\Gal}_* \cong \Gal(\F:\Q^\mathrm{ab})$$
 }

(iv) {\em The group $\bar{N}^\Gal_* $ acts by the
Galois action on the structure $\pi_0(\mathbf{S})$ of connected components so that
$\bar{T}^{0,\Gal}_* $ fixes each component,  and its quotient
$ \hat{\Z}^\times$
acts on $\pi_0(\mathbf{S})$ transitively.}

{\bf Proof.}  Since
$$\Gal(\F: \Q)= \varprojlim\Gal(\F_{\bar{\tau}}:\Q),$$
the isomorphisms (i) follow from the isomorphism statements in 
\ref{318} and \ref{Cor321}.  For the same reason (ii) follows from \ref{N/T}.  (iii) follows from 
\ref{320} and \ref{(b)}.

(iv) follows from the description of  components in \ref{defApprox} in terms of determinants along with the statement \ref{TT1}(i). 

$\Box$
\epk
\section{The structure on $\HH.$ }
\bpk
Here we aim to describe the structure on $\tilde{\HH},$ the projective limit of $\Gamma(N)\backslash \HH.$ As pointed out by the bijection (\ref{Sapprox})  $\tilde{\HH}$ is effectively equal to
$\Ss_\approx$ (see \ref{defApprox}--\ref{S0+}), the quotient of $\Ss$ by the action of $\Delta(\hat{\Z}).$

The structure on the set $\Ss_\approx$  which is of interest to us is exactly the structure which is induced (in model-theoretic sense) from $\mathbf{S}^\mathrm{Pure}.$ We call the structure
 $\mathbf{S}^\mathrm{Pure}_\approx.$ The exact meaning of the definition is that we aim to identify the structure $\mathbf{S}^\mathrm{Pure}_\approx$ such that the quotient map $\Ss\twoheadrightarrow \Ss_\approx$ induces  
the epimorphism
\be\label{defSapprox} \Aut\, \mathbf{S}^\mathrm{Pure}\twoheadrightarrow \Aut\, \mathbf{S}^\mathrm{Pure}_\approx.\ee

 Even more specifically we are interested in the substructure $\mathbf{S}^\mathrm{Pure}_\approx(\mathrm{CM})$ on its special points.

 Similarly we define  $\mathbf{S}^\mathrm{Full}_\approx$ as the structure on $\mathbf{S}_\approx$ induced from
$ \mathbf{S}^\mathrm{Full}$ with respective epimorphism of the automorphism groups
\be\label{defSapproxFull} \Aut\, \mathbf{S}^\mathrm{Full}\twoheadrightarrow \Aut\, \mathbf{S}^\mathrm{Full}_\approx.\ee

\epk
\bpk \label{Autapprox} {\bf Proposition.}  {$\Aut\, \mathbf{S}^\mathrm{Pure}(\mathrm{CM})$ acts faithfully on $\mathbf{S}^\mathrm{Pure}_\approx(\mathrm{CM}).$ That is 
$$\Aut\, \mathbf{S}^\mathrm{Pure}_\approx(\mathrm{CM})=\Aut\, \mathbf{S}^\mathrm{Pure}(\mathrm{CM}).$$}

{\bf Proof.} By definition
\be\label{gD0} [\tau,g\Delta]^\sigma=[\tau,r_\tau(\sigma)\cdot g\Delta].\ee

We need to show that for $\sigma\in \Aut\, \mathbf{S}^\mathrm{Pure}$ the condition
\be\label{gD} [\tau,g\Delta]^\sigma=[\tau,g\Delta]\ee
for all $\tau\in \HH$ quadratic and all $g\in \GL_2 ,$ implies that
$\sigma$ is  the identity (recall the notation in \ref{Del}).

We assume first that in (\ref{gD}) $\sigma$ acts trivially on $\Ss_\approx.$ $g=1$ and
 that $\tau=\sqrt{-n}.$  
 
 By definition
$$[\tau,1]^\sigma=[\tau,r], \mbox{ for some }r=
\left(\begin{array}{ll}
\pm a\ \pm nb\\ -b\ \ \ \ a
\end{array}\right)$$
It follows $r\in \Delta$ and hence $b=0$ and $a=1,$ that is $r=\rrr_{-1}$ or the identity. We assume, towards a contradiction, that $\sigma$ acts non-trivially on $S_{\tau}$ and thus $r=\rrr_{-1}$.

Now for the more general case we have with the same $r$
$$[\tau,g]^\sigma=[\tau,r\cdot g]=[\tau,q\cdot r\cdot g]$$
for any $q\in G_{\tau}$
and thus (\ref{gD}) implies 
$$g\inv\cdot q\rrr_{-1}\cdot g\in \Delta,\mbox{ for some }q\in G_{\tau} $$
But $\det q\in \Q_+$ and  $\det g\inv\cdot q\rrr_{-1}\cdot g=\det q\rrr_{-1}$ is a negative rational number. In $\Delta$ this is only satisfied for  $g\inv\cdot q\rrr_{-1}\cdot g=\rrr_{-1}.$

Thus $q\rrr_{-1}=\left(\begin{array}{ll}
-a\ - nb\\ -b\ \ \ \ a
\end{array}\right)\in  \GL_2(\Q).$ This is an involution only when $b=0$ and $a=\pm 1.$ Now the extra condition 
$g\inv\cdot q\rrr_{-1}\cdot g=\rrr_{-1}$ for all $g$ gives the final contradiction. $\Box$

\epk

\bpk \label{specpointsapprox} {\bf The  structure $\mathbf{S}^\mathrm{Pure}_\approx$. } We describe explicitly the basic relations of the structure.

A. {\em Curves  on $\mathrm{S}_\approx$
 and
  on $\mathrm{S}_{\approx_{K(N)}}$  induced by the action of $\GL_2 .$ }
 
 Let $h$ be an element of $\mathrm{GL}_2 $ and let, 
 $$C_{h}=\{ \la x,h*x\ra : x\in \Ss\}\subset \Ss^2, $$
and, for $K=K(N),$  
$$C_{h,K}=\{ \la x_K,(h*x)_K\ra : x\in\Ss,\ x_K=xK\in \Ss_K\}\subset \Ss^2_K. $$
The latter is the image of $C_h$ in $\Ss_K^2.$ Applying 
 $\approx$ to $\Ss$ and $\Ss_K$ 
 we get  curves $$C_{h/\approx} \subset  \Ss^2_\approx\mbox{ and }C_{h/\approx_K} \subset  \Ss^2_{\approx_K}.$$ Recall that   $\Ss_{\approx_K}$ is the curve isomorphic to the irreducible curve $\Y(N)$.
$C_{h/\approx_K}$  is in general reducible and splits into irreducible components  $C^i_{h/\approx_K}.$ 
We remark that, by construction: 

(a) {\em The curve $C_{h/\approx_K}\subset \Ss^2_{\approx_K}$ and its splitting into irreducible components  $C^i_{h/\approx_K}$ are invariant under the action of $\Aut\, \C.$

Respectively, the  curve  $C_{h/\approx}\subset \Ss^2_\approx$ and its splitting into irreducible components  $C^\mu_{h/\approx}$ 
are invariant under the action of $\Aut\, \C.$
} 

\medskip

Let $\tilde{\G}=\Delta(\Q_+)\cdot \SL_2 ,$ a subgroup of $\GL_2 .$

Below  we narrow our analysis to the case when $h\in \tilde{\G}.$ Note that  
in this case the components $\Ss^\mu$ of $\Ss$ remain invariant under 
   the map 
  $x\mapsto h*x$
  and hence  the curve $C_h$ splits 
    $$C_h=\bigcup_{\mu\in \hat{\Z}^\times} C_h^\mu; \ \ C_h^\mu= C_h\cap (\Ss^\mu\times \Ss^\mu),$$
and after factoring by $\approx$
 $$C_{h/\approx}=\bigcup_{\mu\in \hat{\Z}^\times} C_{h/\approx}^\mu;\ \ C_{h/\approx}^\mu\subset \Ss_\approx\times \Ss_\approx$$
 where the canonical projections
 $$   C_{h}^\mu \to   C_{h/\approx}^\mu$$
 are bijections. It follows that 
 
(b) {\em The curves  $C_{h/\approx}^\mu\subset \Ss_\approx\times \Ss_\approx$
 are the graphs of the bijective maps $$h^\mu: \Ss_\approx\to  \Ss_\approx.$$ Any $\sigma\in \Aut\, \C$ acts on $\bigcup_\mu h^\mu$ so that
 $(h^\mu)^\sigma=h^\lambda,$ some $\lambda\in \hat{\Z}^\times,$
 and thus the
 union $\bigcup_\mu h^\mu$ is invariant under the action of $\sigma.$
   } 
 
Note  that the natural action of $\rrr_\lambda\in \Delta$ on $\Ss_\approx$ induced from $\Ss$ fixes  $\Ss_\approx$ point-wise. However, in accordance with the  above, 
$$ (\rrr_\lambda\inv\cdot h\cdot \rrr_\lambda)^\mu= h^{\mu\lambda}.$$
 In other words, when starting with a given $h\in \tilde{\G}$  the invariant set of transformations $h^\mu:  \Ss_\approx\to \Ss_\approx$  corresponds to the conjugacy class 
 $$h^\Delta=\{ \rrr_\lambda\inv\cdot h\cdot \rrr_\lambda: \ \rrr_\lambda\in \Delta\}.$$
 $h$ is one of its indistinguishable elements.

The   conjugacy class $h^\Delta$ corresponds to the curve $C_h$ and 
$h$ corresponds to one of the irreducible components of $C_h.$ 
   
 To sum up:
\epk
 \bpk\label{4.3(c)} {\bf Proposition.}
  {\em $\tilde{\G}$ as the group acting 
 on $\Ss_\approx$ can be represented by the family of curves $C_{h/\approx},$ $h\in \tilde{\G},$ with individual elements of $\tilde{\G}$ represented by irreducible components $C_{h/\approx}^\mu$ of $C_{h/\approx}.$ The $C_{h/\approx}^\mu$ are graphs of the action on $\Ss_\approx$ of the respective elements of   $\tilde{\G}.$} 
 
 {\em Any element of $s\in \Ss_\approx$ can be represented canonically 
 by an element\linebreak $[\tau,a]\in \Ss,$ \be\label{srep} s=[\tau,a\Delta],\  
 \tau\in \HH, \ 
\det a\in \Q_+.\ee
The choice of 
such an  $a$ for $s$ is unique.

The action of $g\in \tilde{\G}$ on $s:$
\be\label{sact} g*[\tau,a\Delta]=[\tau,ag\inv \Delta].\ee
 }

{\bf Proof.} First note that any coset $a\Delta\subset \GL_2(\A_f),$ $a\in  \GL_2(\A_f),$
contains exactly one element of determinant in $\Q_+.$ The uniquiness of the representation  (\ref{srep}) follows.

The action (\ref{sact}) is well-defined because of the uniquiness of $a$ and the fact that $\det g\in \Q_+.$

Finally, note that the action (\ref{sact}) is isomorphic to the action of $g$ on the component $\Ss^1$ of $\Ss,$ which according to \ref{specpointsapprox}(b) represents the action of $g$ on $\Ss_\approx.$
$\Box$

\medskip

{\bf Remark} ({\em Diagonal matrices and their action})

Let $$\rrr_{a,b}:= \left( \begin{array}{ll}
a\ 0\\
0\ b
\end{array}\right), \mbox{ for }a,b\in \A_f.$$


Since the $\rrr_{a,b}$ commute with $\rrr_\mu,$ the  action by $\rrr_{a,b}$
on $\Ss$ passes to an action on $\Ss_\approx,$ 
$$x\mapsto \rrr_{a,b}\cdot x$$ and the action is invariant under automorphisms. 

It is immediate by construction that $\rrr_\mu\cdot x=x,$ for $x\in \Ss_\approx,$  $\mu\in \hat{\Z}^\times.$
 Hence 
\be\label{dmatr} \rrr'_\mu\cdot x= \rrr_{\mu\inv, \mu}\cdot x. \ee
 And note that $\rrr_{\mu\inv, \mu}\in \SL_2(\A_f).$
 
 \medskip  

{\bf Warning.} Although both  $\tilde{\G}=\Delta(\Q_+)\cdot \SL_2(\A_f)$
and $\Delta'(\A_f)$ act on $\Ss_\approx,$ the action does not extend to the group $\GL_2(\A_f)$ generated by the two subgroups.

\epk  


\bpk\label{C} {\bf The structure $(\Ss_\approx,  \tilde{\G}).$}
The action of $\tilde{\G}$ on $\Ss_\approx$ induces 
the action of $\mathrm{SL}_2(\hat{\Z})$ on $\mathrm{S}_{\approx_{K(N)}},$ where $K(N)=\Delta\cdot \tilde{\Gamma}(N)=  \tilde{\Gamma}(N)\cdot \Delta.$

Since a definition of an individual element $h^\mu$ of $\tilde{\G}$ requires the definition of a component $\Ss^\mu$ of $\Ss,$  the action of elements of $\mathrm{SL}_2(\hat{\Z})$ on $\mathrm{S}_{\approx_{K(N)}}$
is defined over $\Q^\mathrm{ab}$. 

Note that for each $N$  the subgroup $\tilde{\Gamma}(N)\subset \tilde{\G}$ can be distinguished as the set of those elements which act as identity on $\Ss_{\approx_{K(N)}}.$

\medskip
 
 {\em Special points.} 
It is immediate from \ref{4.3(c)} and \ref{specpoints} that 
  the fixed points $s$ of  elements of   $g\in \tilde{\G}$ in $\Ss_\approx$ can be represented as $s=[\tau,a\cdot \Delta],$ $\det a\in \Q_+,$
  for $\tau$ quadratic such that 
$g*[\tau,a]=[\tau,a].$ 


Respectively, the special points of  $\mathrm{S}_{\approx_{K}}$ are of the form $[\tau,a\cdot\Delta\cdot K],$ for $\tau\in \HH$ quadratic.

\medskip
{\em 4-point relations between special points} (induced by the equivalence relation $E,$ see \ref{definitionS}(iii)). 

For each quadratic $\tau_1, \tau_2$ 
define the 4-ary relation $R_{\tau_1,\tau_2}(s_1,s_2,t_1,t_2)$ on $\Ss_\approx:$
$$R_{\tau_1,\tau_2}(s_1,s_2,t_1,t_2)\Leftrightarrow $$
$$ \Leftrightarrow \bigwedge_{i=1,2}s_i=[\tau_i,a_i]_\approx \ \& \ t_i=[\tau_i,b_i]_\approx\ \& \
\exists r_1\in N_{\tau_1}, r_2\in N_{\tau_2}\ $$ $$ \det r_1=\det r_2\in \Delta \ \& \
\bigwedge_{i=1,2} b_i=r_i\cdot a_i\cdot \det r_i\inv
$$
where we assume that the representation of the $s_i$ and $t_i$ is canonical and $\det r_i$ takes its values in $\Delta$ via the isomorphism $\Delta\cong \hat{\Z}^\times.$ Note that we may always assume
that $\tau_1=\sqrt{-m_1}$ and $\tau_2=\sqrt{-m_2}$ for positive square-free integers.

Now set 
$$R(s_1,s_2,t_1,t_2)\equiv \bigvee R_{\sqrt{-m_1},\sqrt{-m_2}}(s_1,s_2,t_1,t_2)$$
where $m_1,m_2$ run through square-free positive integers.
\medskip

(a) {\em Relation  $R$ is invariant under $\Aut\, \mathbf{S}^\mathrm{Pure}_\approx.$   } 
Indeed, by \ref{Autapprox} any $\sigma\in \Aut\, \mathbf{S}^\mathrm{Pure}_\approx$ is induced by the unique $\sigma\in \Aut\, \mathbf{S}^\mathrm{Pure}.$ The latter acts on special points by
$$[\tau_i,a_i]\mapsto [\tau_i,r_{\tau_i}(\sigma)\cdot a_i]; \ \
 [\tau_i,b_i]\mapsto [\tau_i,r_{\tau_i}(\sigma)\cdot b_i],$$
where $r_{\tau_i}(\sigma)\in N_{\tau_i},$ $\det r_{\tau_1}(\sigma)=\det r_{\tau_2}(\sigma).$  Respectively, in terms of  $\mathbf{S}^\mathrm{Pure}_\approx$ the result of the application of $\sigma$ are
 $[\tau_i,a'_i]_\approx$ and $[\tau_i,b'_i]_\approx$ where
 $$a'_i=r_{\tau_i}(\sigma)\cdot a_i\cdot \det r_{\tau_i}(\sigma)\inv; \ \
 b'_i=r_{\tau_i}(\sigma)\cdot b_i\cdot \det r_{\tau_i}(\sigma)\inv$$
 Assuming $R_{\tau_1,\tau_2}(s_1,s_2,t_1,t_2),$ for $s_i=[\tau_i,a_i]_\approx\ \& \ t_i=[\tau_i,b_i]_\approx$ we get by definition
 $$b'_i=(r_{\tau_i}(\sigma)\cdot r_i\cdot r_{\tau_i}(\sigma)\inv)\cdot a'_i\cdot \det (r_{\tau_i}(\sigma)\cdot r_i\cdot r_{\tau_i}(\sigma)\inv)\inv $$
 and $r_{\tau_i}(\sigma)\cdot r_i\cdot r_{\tau_i}(\sigma)\in N_{\tau_i}.$
 Thus $R_{\tau_1,\tau_2}(s_1^\sigma,s_2^\sigma,t_1^\sigma,t_2^\sigma)$ holds.
 
 (b) {\em  For any special $s_1=[\tau_1,a_1]_\approx,\ s_2=[\tau_2,a_2]_\approx$ and any $\sigma\in \Aut\, \mathbf{S}^\mathrm{Pure}$ it holds
$R(s_1,s_2, s_1^\sigma,  s_2^\sigma)$, where $s^\sigma$ is the  induced action.}

Indeed, we have
$$\sigma: [\tau_i,a_i]_\approx\mapsto [\tau_i,r_{\tau_i}(\sigma)\cdot a_1]_\approx$$
and in terms of canonical representatives
 $$\sigma: [\tau_i,a_i]\mapsto [\tau_i,r_{\tau_i}(\sigma)\cdot a_1\cdot \det r_{\tau_i}(\sigma)\inv]$$
 as required by $R_{\tau_1,\tau_2}.$ $R(s_1, s_1^\sigma, s_2, s_2^\sigma)$ follows.

(c) {\em For any special $s_1=[\tau_1,a_1]_\approx,\ s_2=[\tau_2,a_2]_\approx,$

$R(s_1,s_2, t_1,  t_2)$ iff there is  $\sigma\in \Aut\, \mathbf{S}^\mathrm{Pure}$ such that $t_1=s_1^\sigma$ and $t_2=s_2^\sigma.$}

Proof. (b) proves one direction. In the opposite direction use the fact that by definition 
$$t_i=[\tau_i, r_i\cdot a_i\cdot \det r_i\inv]_\approx=[\tau_i, r_i\cdot a_i]_\approx$$
where $\det r_1=\det r_2 \in \Delta.$ By \ref{Cor316} there is $\sigma\in \Aut\, \mathbf{S}^\mathrm{Pure}$ such that
$r_i=r_{\tau_i}(\sigma)$ modulo $\Q_+.$ Using the Hasse principle as in the proof of \ref{TT1}(ii) we can make this an equality. This proves the claim.


\epk


\bpk \label{Hpure}
 It is  useful to think of the structure $\mathbf{S}^\mathrm{Pure}_\approx$ as the structure on the universe $\tilde{\HH}$ which by  (\ref{Sapprox}) can be identified with $\Ss_\approx.$ The group $\tilde{\G}$ acts on $\tilde{\HH}$ but the action is defined up to  inner automorphisms  as described in \ref{specpointsapprox}A(c). The other relations on $\tilde{\HH}$ described in \ref{specpointsapprox}, B - E, determine a structure on $\tilde{\HH}$ which we call $\tilde{\HH}^\mathrm{Pure}.$

\medskip

{\bf Definition.} $\tilde{\HH}^\mathrm{Pure}$ is the structure with two universes (sorts) $\Ss_\approx$ 
and $\tilde{\G}$ (the set $\Ss_\approx$ will also be called  $\tilde{\HH}$); 

- $\tilde{\G}$ is a group isomorphic to the product $\Delta(\Q_+)\cdot \SL_2(\A_f);$ it  acts on $\Ss_\approx$ by the action induced from $\mathbf{S}^\mathrm{Pure};$

- for each $g\in \tilde{\G}$ the set $$g^\Delta=\{ \varphi_\mu (\mathbf{g}):\ \mu\in \hat{\Z}^\times\},\mbox{ for } \mathbf{g}\in \Delta(\Q_+)\cdot \SL_2(\A_f),\mbox{ such that }g=\varphi_1(\mathbf{g}),$$ is definable (that is invariant under automorphisms);


- the subgroups $\tilde{\Gamma}(N)\subset  \tilde{\G}$ are definable;

- for any special $s_1,s_2\in \Ss_\approx$ the set $$\mathrm{tp}(\la s_1,s_2\ra):=\{ \la t_1,t_2\ra\in \Ss_\approx^2: R(s_1,s_2,t_1,t_2)\}$$ is definable (the model-theoretic notation for the {\em type} of the pair of points).

\epk

\bpk \label{Gapprox} {\bf Theorem}. {\em The structure $\mathbf{S}^\mathrm{Pure}_\approx$ 
 is bi-interpretable with the structure $\tilde{\HH}^\mathrm{Pure},$ and $\mathbf{S}^\mathrm{Pure}_\approx(\mathrm{CM})$ is bi-interpretable with $\tilde{\HH}^\mathrm{Pure}(\mathrm{CM}).$
 
  }

{\bf Proof.} By definition, see (\ref{defSapprox}), we need to prove 
that  any automorphism of $\mathbf{S}^\mathrm{Pure}$ induces an automorphism of $\tilde{\HH}^\mathrm{Pure}$ and, conversely, any automorphism of  $\tilde{\HH}^\mathrm{Pure}$
 can be lifted to an automorphism of 
$\mathbf{S}^\mathrm{Pure}.$ The former is proved in \ref{specpointsapprox}. So we need to show that any  
$\sigma\in \Aut\,\tilde{\HH}^\mathrm{Pure}$ can be lifted to an automorphism $\rho$ of 
$\mathbf{S}^\mathrm{Pure}.$ 

Claim. {\em Any 
$\sigma\in \Aut\,\tilde{\HH}^\mathrm{Pure}(\mathrm{CM})$ can be lifted to an automorphism $\rho$ of 
$\mathbf{S}^\mathrm{Pure}(\mathrm{CM}).$ }

Proof of Claim. Consider special points
$s_1,\ldots,s_n,\ldots\in \Ss_\approx$ with representatives \linebreak
$[\tau_1, a_1],\ldots,[\tau_1, a_n],\ldots\in \Ss,$
$$\sigma: s_i\mapsto t_i,\ \ t_i=[\tau_i,b_i]_\approx, \mbox{ canonical representatives}.$$

Since $\sigma$ preserves $R,$ we get 
$$\vDash  R(s_1,s_i,t_1,t_i), \mbox{ for }i=1,...,n,...$$
and hence
$$b_i=r_i\cdot a_i\cdot \det r_i\inv,\ r_i\in N_{\tau_i}, \ \det r_1=\det r_i, \mbox{ for }i=1,...,n...$$
By  \ref{Cor316} there is a $\rho\in  \Aut\, \mathbf{S}^\mathrm{Pure}(\mathrm{CM})$ such that $r_i=r_{\tau_i}(\rho),$ $i=1,...,n...$,  and hence
$$s_i^\sigma=t_i=[\tau_i,b_i \Delta]=[\tau_i,a_i\Delta]^\rho=s_i^\rho,\ \ i=1,...,n... $$
Claim proved.
 
Note that the last argument in the proof of the claim also shows that $\bigwedge_i R(s_1,s_i,t_1,t_i)$ implies
the existence of $\sigma: s_i\mapsto t_i.$ 

\medskip

Now let $\sigma\in \Aut\, \tilde{\HH}^\mathrm{Pure}$ and $\sigma^\mathrm{CM}$ be its restriction to  $\tilde{\HH}(\mathrm{CM}).$ By the claim there is $\rho^\mathrm{CM}\in \Aut\, \mathbf{S}^\mathrm{Pure}(\mathrm{CM})$ whose restriction to $\tilde{\HH}(\mathrm{CM})$ is $\sigma^\mathrm{CM}.$ 

Since special points are Zariski dense in the components $\Ss^\mu_{K(N)},$ $\rho^\mathrm{CM}$ determines the unique bijective map  
$\Ss^\mu\to \Ss^{\rho^\mathrm{CM}(\mu)}$ on the set of components of $\Ss.$ More precisely, 
$$\rho^\mathrm{CM}(\mu)=\mu\cdot \lambda_{\rho}$$
where $\lambda_{\rho}=\det r_\tau(\rho^\mathrm{CM})\in \hat{\Z}^\times$ is obtained from the property (\ref{rsig}) of automorphisms acting on CM-points.

Consider an arbitrary $s\in \Ss_\approx.$ It follows from the definition of $\approx$ that for each component  $\Ss^\mu$ there is a unique 
$x\in \Ss^\mu$ such that $x_\approx=s.$ Call this element $x(s,\Ss^\mu).$ 
By construction any element of $\Ss$ has the form $x(s,\Ss^\mu)$ for some
$s\in \Ss_\approx$ and a component $\Ss^\mu$ of $\Ss.$

Consider the element $\ttt\in \GL_2^+(\Q)$ and the set $\ttt^\Delta$  definable in $\tilde{\HH}^\mathrm{Pure}$ by definition. 
 Elements $\ttt^\mu$ of $\ttt^\Delta$ can be identified with curves $C^\mu_{\ttt/\approx}$ on $\tilde{\HH}\times \tilde{\HH}$ (see \ref{specpointsapprox}A(b)) and there is a one-to-one correspondence
between elements $\ttt^\mu$ of  $\ttt^\Delta$ and components $\Ss^\mu$ of $\Ss$   given by
$$C^\mu_{\ttt}=C_{\ttt}\cap (\Ss^\mu\times \Ss^\mu)$$
(see \ref{specpointsapprox}A).

Clearly, $$\rho^\mathrm{CM}: \ttt^\mu\mapsto \ttt^{ \rho^\mathrm{CM}(\mu)}=\sigma(\ttt^\mu)=\ttt^{ \sigma(\mu)},\ \ \ttt^{ \rho^\mathrm{CM}(\mu)}=\ttt^{\mu\cdot \lambda_\rho}=\ttt^{\mu\cdot \lambda_\sigma}.$$
for $\sigma$ and $\rho^\mathrm{CM}$  as above. This gives us an extension of the action of $\sigma\in \Aut\, \tilde{\HH}^\mathrm{Pure}$ to components of $\Ss,$ $$\sigma: \Ss^\mu \to \Ss^{\mu\cdot \lambda_\sigma}.$$

Define the map $$\rho: \Ss\to \Ss; \ \ x(s,\Ss^\mu)^\rho:= x(s^\sigma, \Ss^{\mu\cdot \lambda_\sigma}).$$
This is the lift of $\sigma\in \Aut\, \tilde{\HH}^\mathrm{Pure}$ to a transformation  $\rho$ of $\Ss$ and by definition the restriction of $\rho$ to the CM-substructure is $\rho^\mathrm{CM}$   above.

By construction $\rho$ respects the splitting of $\Ss$ into components. 

In order to prove that $\rho\in \Aut\, \mathbf{S}^\mathrm{Pure}$
it remains 
to prove that $\rho$ respects the action by elements of $\GL_2(\A_f).$

Consider $h\in \Delta'(\Q_+)\cdot \SL_2(\A_f).$
By definition \ref{specpointsapprox}A(b)) via the action on $\Ss^\mu$  $h$ induces on $\tilde{\HH}$ the transformation
$h^\mu.$ This means that
$$h\cdot x(s,\Ss^\mu)=x(h^\mu\cdot s,\Ss^\mu)$$
and thus
$$(h\cdot x(s,\Ss^\mu))^\rho=x(h^\mu\cdot s,\Ss^\mu)^\rho=
x((h^{\mu}\cdot s)^\sigma, \Ss^{\mu\cdot \lambda_\sigma})=x(h^{\mu\cdot \lambda_\sigma}\cdot s^\sigma, \Ss^{\mu\cdot \lambda_\sigma})=$$ $$=h\cdot x(s^\sigma, \Ss^{\mu\cdot \lambda_\sigma})=
h\cdot x(s,\Ss^\mu)^\rho$$
which proves that any such $h$ is invariant under $\rho.$

Next consider the action of an element $d$ of $\Delta(\A_f)$ on 
$x(s,\Ss^\mu).$ Since any $a\in \A_f$ can be represented as $a=q\cdot \nu,$ for $q\in \Q_+$ and $\nu\in \hat{\Z}^\times,$
$d$ can be represented as 
$d=\rrr_q\cdot \rrr_\nu.$ Thus, the action by $d$ is invariant under $\rho$ if and only if the action by $\rrr_\nu$ is. 

We have by construction
$$\rrr_\nu\cdot x(s,\Ss^\mu)=x(s,\Ss^{\mu\cdot\nu})$$

Thus
$$(\rrr_\nu\cdot x(s,\Ss^\mu))^\rho= x(s,\Ss^{\mu\cdot\nu}))^\rho=x(s^\sigma,\Ss^{\mu\cdot\nu\cdot \lambda_\sigma})=\rrr_\nu\cdot x(s^\sigma,\Ss^{\mu\cdot \lambda_\sigma})=\rrr_\nu\cdot x(s,\Ss^{\mu})^\rho$$
which proves that the actions by $\rrr_\nu$ and by $d$ on $\Ss$ are invariant under $\rho.$

Since any element of $\GL_2(\A_f)$ can be represented in the form
$d\cdot h$ where $d\in \Delta(\A_f)$ and $h\in \SL_2(\A_f)$
we proved that $\rho$ satisfies (\ref{autg}) and thus $\rho\in \Aut\, \mathbf{S}^\mathrm{Pure}.$ $\rho$ is the lifting of $\sigma\in \Aut\, \tilde{\HH}^\mathrm{Pure}$ to  $\Aut\, \mathbf{S}^\mathrm{Pure}.$

This completes the proof of the theorem. $\Box$

\epk

\bpk \label{RemCor} {\bf Remark/Corollary.} 
Consider the embedding 
$$(\HH,\G)\subset (\tilde{\HH},\tilde{\G}),$$
where $\G\cong\GL_2(\Q)\cap \tilde{\G}=\GL_2^+(\Q).$

1. 
Clearly, under this embedding special points in $\HH$ are just quadratic imaginary points. 

The action of  individual elements $\rrr'_q\in \Delta'(\Q_+),$  $\rrr_q\in \Delta(\Q_+)$ and the scalar matrices $q\cdot \mathbf{I}$ are definable. 

The scalar matrices act trivially, as identity.

2. More generally, the action of $h\in \GL_2^+(\Q)$ depends on the choice of an isomorphism   
$\phi_\mu: \G\to \GL_2^+(\Q),$ for   $\mu\in \hat{\Z}^\times,$ see \ref{specpointsapprox} A(b). However, on $\HH$ only
$h^\mu$ for $\mu=\pm 1$ corresponds to a non-empty curve. Thus $h$ defines (in the natural $L_{\omega_1,\omega}$-language) the two-component curve $C_h\cup C_h^{-1},$ where  $C_h^{-1}$ is the graph
of $\rrr_{-1}\cdot h\cdot \rrr_{-1}$
   and the choice of the component depends on $\phi_\mu.$ (Note that if
   $h=\left( \begin{array}{ll}
a\ b\\ c\ d
\end{array}\right)$   then $\rrr_{-1}\cdot h\cdot \rrr_{-1}=\left( \begin{array}{ll}
 a\ -b\\ -c\ \ d
\end{array}\right).$)

\medskip

3. The 4-point relation $R(s_1,s_2,s'_1,s'_2)$ on $\HH$ for
$s_i=p_i+q_i\sqrt{-m_i},$ $s'_i=p'_i+q'_i\sqrt{-m_i},$ $p_i,q_i,p'_i,q'_i\in \Q,$ $q_i,q'_i\in \Q_+,$ holds if and only if $$q_1=q'_1 \ \& \ q_2=q'_2\ \& \ ((p_1=p'_1 \ \& \ p_2=p'_2) \vee (p_1=-p'_1 \ \& \ p_2=-p'_2))$$
   
 Thus in $\HH:$
   $$R(s_1,s_2,s'_1,s'_2)\ \Leftrightarrow\ \la s'_1,s'_2\ra =\la s_1,s_2\ra \ \vee \ \la s'_1,s'_2\ra =\la -s^c_1,-s^c_2\ra$$ 
where $-s^c$ is the application of complex conjugation to $-s.$
\medskip

4. Respectively, on $\Y(N)=\Gamma(N)\backslash \HH=\tilde{\Gamma}(N)\backslash \tilde{\HH}$ the indiscernible curves $C^\mu_{h,K(N)}$
are conjugated by elements of the group $\Delta/\Delta(N),$ or equivalently, by $\Gal(\Q(\zeta_N):\Q).$ The union  $\bigcup_\mu C^\mu_{h,K(N)}$ is a (reducible) algebraic curve over $\Q.$ 

The $j_N$-image of the 4-point relations $R(s_1,s_2,s'_1,s'_2)$ is a relation on CM-points  of $\Y(N),$  call it  $R_N.$
 
 The respective set of pairs (called type of $\la j_N(s_1),j_N(s_2)\ra$ in Model Theory), 
$$\mathrm{tp}_N(\la j_N(s_1),j_N(s_2)\ra):=\{ \la y_1,y_2\ra\in \Y(N)\times \Y(N): R_N(j_N(s_1),j_N(s_2),y_1,y_2)\}$$
is, by \ref{specpointsapprox} E(c), a $\Gal_\Q$-orbit of the pair, a finite Zariski set. 

\medskip 

5. Given $g\in \tilde{\G}$ the graph of $g,$ $\mathsf{graph}\, g\subset \tilde{\HH}\times \tilde{\HH},$ has a point in $\HH\times \HH$ if and only if $g\in \GL_2^+(\Q).$ 

Indeed, suppose $\la s,s'\ra\in (\HH\times \HH)\cap \mathsf{graph}\, g.$
We can represent $s=[\tau, \Delta]$  and 
$s'=[\tau',\Delta]$ $\tau,\tau' \in \HH.$ Then by assumption
$[\tau',\Delta]=[\tau,g\inv\Delta],$ which implies 
$\tau'=q\cdot \tau$ and 
$[q\cdot\tau,q\cdot g\inv\Delta]=[q\cdot\tau,\Delta],$  for some $q\in \GL_2^+(\Q).$  That is
$q\cdot g\inv \in \Delta$ which can only happen if $q\cdot g=1,$ since
$\det q\cdot g\in \Q_+.$

\medskip

\epk

\bpk {\bf The structure $\tilde{\HH}^\mathrm{Full}.$} The determination of this structure in terms of concrete relations and operations as  in \ref{Hpure} for $\tilde{\HH}^\mathrm{Pure}$
would be rather cumbersome for the context of the current paper. Instead we recall that the respective structure $\mathbf{S}^\mathrm{Full}$ is the structure obtained from $\mathbf{S}^\mathrm{Pure}$ by adding some algebraic-geometric relations 
defined over $\Q$ (an {\em expansion} of the structure $\mathbf{S}^\mathrm{Pure}$).
 More technically  the relationship between the smaller and the bigger structure is described by the close relationship between their automorphism groups  in section 3. 
 
 \medskip
 
 {\bf Definition.} $\tilde{\HH}^\mathrm{Full}$ is the expansion of $\tilde{\HH}^\mathrm{Pure}$ 
 defined over $\Q,$ such that
  $$\Aut\, \tilde{\HH}^\mathrm{Full}=\Aut\, \mathbf{S}^\mathrm{Full}.$$
  
\epk
\bpk {\bf Corollary.} {\em $\Aut\, \C$ acts on $\tilde{\G}$ in $\tilde{\HH}^\mathrm{Full}$ as 
$\Gal(\Q^\mathrm{ab}:\Q)\cong \hat{\Z}^\times:$}

for $\mu\in \hat{\Z}^\times$ and $h=\left(\begin{array}{ll}
a \ b\\ c \ d
\end{array}\right)\in \tilde{\G},$ \ $h^\mu=\left(\begin{array}{ll}
a \ \ \ \mu b\\ \mu\inv c \ \ d
\end{array}\right).$

{\bf Proof.} See \ref{specpointsapprox}, A(b). 
\epk

\end{document}